\documentclass[11pt]{amsart}
\usepackage{amssymb}

\newtheorem{lemma}{Lemma}[section]

\newtheorem{theorem}[lemma]{Theorem}
\newtheorem{remark}[lemma]{Remark}
\newtheorem{proposition}[lemma]{Proposition}
\newtheorem{corollary}[lemma]{Corollary}
\newtheorem{definition}[lemma]{Definition}

\begin{document} 
\title{Radial index and Poincar\'e-Hopf index of 1-forms on semi-analytic sets}
\author{Nicolas Dutertre}
\address{Universit\'e de Provence, Centre de Math\'ematiques et Informatique,
39 rue Joliot-Curie,
13453 Marseille Cedex 13, France.}
\email{dutertre@cmi.univ-mrs.fr}

\thanks{
 Supported by {\em Agence Nationale de la Recherche} (reference ANR-08-JCJC-0118-01)}
\begin{abstract}
 The radial index of a $1$-form on a singular set is a generalization of the classical Poincar\'e-Hopf index.
We consider different classes of closed singular semi-analytic sets in $\mathbb{R}^n$ that contain $0$ in their singular locus
 and we relate the radial index of a $1$-form at $0$ on these sets to Poincar\'e-Hopf indices at $0$ of vector fields defined on
 $\mathbb{R}^n$.

\end{abstract}

\maketitle
\markboth{Nicolas  Dutertre}{Radial index and Poincar\'e-Hopf index of 1-forms on semi-analytic sets}
\section{Introduction}
It is well-known that one can assign to each isolated zero $P$ of a vector field $v$ on a smooth manifold $M$ an index called the
Poincar\'e-Hopf index that we will denote by Ind$_{PH}(v,P,M)$. The Poincar\'e-Hopf theorem says that if $M$ is compact and $v$
admits a finite number of zeros $P_1,\ldots,P_k$ then:
$$\chi(M)=\sum_{i=1}^k \hbox{Ind}_{PH}(v,P_i,M).$$
In [Sc1,Sc2,Sc3] (see also [BrSc]), M-H Schwartz has proved a version of this theorem for a Whitney stratified analytic subvariety of an
analytic manifold $M$ and for a class of vector fields that she called radial vector fields. 
The radial vector fields are defined in terms of two types of tubes around strata. The first tubes are given by the barycentric
subdivision of a triangulation and are called parametric tubes. The second are given by certain geodesic tubular
neighborhoods defined using the ambient metric and are called geodesic tubes.
A radial vector field $v$ is a
continuous vector field on $M$, tangent to the strata of $V$ and exiting from sufficiently small geodesic tubes around the strata of $V$
over closed subsets of the strata that contain the zeros of $v$. Here a point $P$ is a zero of $v$ if it is a zero of $v$
restricted to $X(P)$ where $X(P)$ is the stratum that contains $P$ and the index of $v$ at $P$ is the Poincar\'e-Hopf index of
$v$ restricted to $X(P)$. After the work of M-H Schwartz, several generalizations of the Poincar\'e-Hopf theorem for vector
fields on singular
spaces, together with generalizations of the Poincar\'e-Hopf index, were given (see [ASV], [BLSS], [EG1], [KT], [Si], [SS]). 
The most general version is due to King and Trotman for semi-radial vector fields on radial manifold complexes (Theorem 5.4 in
[KT]). 

Instead of vector fields, one can consider $1$-forms. This is one of the subjects of [Ar], where $1$-forms on manifolds with
boundary are studied. If $M$ is a manifold with boundary $\partial M$ and $\omega$ is a $1$-form, a point in $\partial M$ is a boundary singularity (or a
boundary zero) of $\omega$ if it is a zero of $\omega$ restricted to $\partial M$. To each isolated boundary zero of $\omega$,
Arnol'd assigns an index that he calls the boundary index and proves a Poincar\'e-Hopf theorem for $1$-forms on manifolds with
boundary (see [Ar], p.4). Furthermore, he relates this boundary index to classical Poincar\'e-Hopf indices of vector fields (see
[Ar], p.7). 

In a serie of papers, Ebeling and Gusein-Zade [EG2-6] study $1$-forms on singular analytic spaces. In [EG5], they give a
Poincar\'e-Hopf theorem for a $1$-form on a compact singular analytic set. More precisely, they consider an analytic set $X
\subset \mathbb{R}^N$ equipped with a Whitney stratification and a continuous $1$-form $\omega$ in $\mathbb{R}^N$. A point $P$
in $X$ is a zero (or singular point) of $\omega$ on $X$ if it is a zero of $\omega$ restricted to the stratum that contains
$P$. If $P$ is an isolated zero of $\omega$ on $X$, they define the radial index of $\omega$ at $P$ (Definition p.233 in
[EG5]). Let us
denote it by Ind$_{Rad}(\omega,P,X)$. Then they prove that if $X$ is compact and $\omega$ is a $1$-form on $X$ with a finite
number of zeros $P_1,\ldots,P_k$ then (Theorem 1 in [EG5]) :
$$\chi(X)=\sum_{i=1}^k \hbox{Ind}_{Rad}(\omega,P_i,X).$$ 
It is straightforward to see that the definitions and results of Ebeling and Gusein-Zade extend to the case of closed
subanalytic sets. In this paper, we consider different classes of closed semi-analytic sets in $\mathbb{R}^n$ that contain $0$
in their singular locus and relate the radial index of a $1$-form at $0$ on these sets to classical Poincar\'e-Hopf indices
at $0$ of
vector fields on $\mathbb{R}^n$, like Arnol'd does for manifolds with boundary.

Let us describe the content of the paper. In Section 2, we recall some results about $1$-forms on smooth manifolds. In Section
3, we give a Poincar\'e-Hopf theorem for a class of $1$-forms, called correct, on manifolds with corners (Theorem 3.6). This is not the most
general Poincar\'e-Hopf theorem, as already explained above, but it is enough for our purpose. Moreover, we think that it is worth
stating it in this concrete form. In Section 4, we define the radial index of a $1$-form on a closed subanalytic set. Section 5
is devoted to the study of the radial index on a manifold with corners. Let $(x_1,\ldots,x_n)$ be a coordinate system in
$\mathbb{R}^n$. For $k \in \{1,\ldots,n\}$ and for every $\epsilon=(\epsilon_1,\ldots,\epsilon_k) \in \{0,1\}^k$, let
$\mathbb{R}^n(\epsilon)$ be defined by :
$$\mathbb{R}^n(\epsilon)=\left\{ (x_1,\ldots,x_n)\in \mathbb{R}^n \ \vert \ (-1)^{\epsilon_1} x_1 \ge
0,\ldots,(-1)^{\epsilon_k} x_k \ge 0\right\}.$$ 
We consider a smooth $1$-form $\Omega=a_1dx_1+\cdots+a_ndx_n$ in $\mathbb{R}^n$. Since $\mathbb{R}^n(\epsilon)$ is
semi-algebraic, Ind$_{Rad}(\Omega,0,\mathbb{R}^n(\epsilon))$ is well-defined. In Theorem 5.4, we relate this index to
Poincar\'e-Hopf indices at $0$ of vector fields defined in terms of the $a_i$'s. Section 6 is not related directly to the
radial index but contains results that will be used in Section 7. We consider a smooth vector field $V$ defined in the
neighborhood of the origin in $\mathbb{R}^n$ such that $0$ is an isolated zero of $V$. We assume that $V$ satisfies the
following condition $(P')$: there exist smooth vector fields $V_2,\ldots,V_n$ defined in the neighborhood of $0$ such that
$V_2(x),\ldots,V_n(x)$ span $V(x)^\perp$ whenever $V(x) \not= 0$ and such that $(V(x),V_2(x),\ldots,V_n(x))$ is a direct basis
of $\mathbb{R}^n$. Let $Z$ be another smooth vector field defined in the neighborhood of $0$ and let $\Gamma$ be the
following vector field :
$$\Gamma=\langle V,Z \rangle \frac{\partial}{\partial x_1}+\langle V_2,Z \rangle \frac{\partial}{\partial x_2}+\cdots+\langle V_n,Z
\rangle \frac{\partial}{\partial x_n},$$
where $\langle \ , \ \rangle$ is the euclidian scalar product. The main result of this section is Theorem 6.7, in which we give
an equality between the indices at $0$ of these three vector fields. In Section 7, we consider an analytic function $f
:(\mathbb{R}^n,0) \rightarrow (\mathbb{R},0)$ defined in the neighborhood of $0$ with an isolated critical point at the origin
and a smooth $1$-form $\Omega=a_1dx_1+\cdots+a_ndx_n$. We first assume that $\nabla f$ satisfies Condition $(P')$ above. In
Theorem 7.2 and Theorem 7.6, we relate Ind$_{Rad}(\Omega,0,f^{-1}(0))$, Ind$_{Rad}(\Omega,0,\{f\ge 0\})$ and
Ind$_{Rad}(\Omega,0,\{f \le 0\})$ to Poincar\'e-Hopf indices at $0$ of vector fields defined in terms of $f$ and $\Omega$. Then
we assume that the vector $V(\Omega)=a_1 \frac{\partial}{\partial x_1}+\cdots+a_n \frac{\partial}{\partial x_n}$ dual to
$\Omega$ satisfies Condition $(P')$ and in Theorem 7.10 and Theorem 7.14, we give the versions of Theorem 7.2 and Theorem 7.6
in
this situation. In Section 8, we explain how to compute the radial index of a $1$-form on a semi-analytic curve. More precisely,
let $F=(f_1,\ldots,f_{n-1}) : (\mathbb{R}^n,0) \rightarrow (\mathbb{R}^{n-1},0)$ be an analytic mapping defined in the
neighborhood of the origin such that $F(0)=0$ and $0$ is isolated in 
$\{x \in \mathbb{R}^n \ \vert \ F(x)=0 \hbox{ and } \hbox{rank}[DF(x)]< n-1 \}$. Let $\Omega=a_1dx_1+\cdots+a_ndx_n$ be a smooth
$1$-form and let $g_1,\ldots,g_k : (\mathbb{R}^n,0) \rightarrow (\mathbb{R},0)$ be analytic functions. For every 
$\epsilon=(\epsilon_1,\ldots,\epsilon_k) \in \{0,1\}^k$, let $\mathcal{C}(\epsilon)$ be the semi-analytic curve
defined by :
$$\mathcal{C}(\epsilon)=F^{-1}(0) \cap \{(-1)^{\epsilon_1}g_1 \ge 0,\ldots,(-1)^{\epsilon_k}g_k \ge 0\}.$$
In Theorem 8.4, Corollary 8.5 and Theorem 8.6, we express the indices Ind$_{Rad}(\Omega,0,F^{-1}(0))$ and
Ind$_{Rad}(\Omega,0,\mathcal{C}(\epsilon))$ in terms of Poincar\'e-Hopf indices at $0$ of vector fields defined in function of
$\Omega$, $F$ and the $g_i$'s. 

When the vector fields that appear in our results have an algebraically zero at $0$, we can apply the
Eisenbud-Levine-Khimshiashvili formula ([EL], [Kh]) and obtain algebraic formulas for the radial index of a $1$-form. One should mention that
this aspect of our work is related to the work of several authors on algebraic formulas for the GSV-index, which is another
generalization of the Poincar\'e-Hopf defined in [GSV] (see [EG2], [EG3], [GGM], [GM1], [GM2], [Kl]). 

Some explicit computations are given to illustrate our formulas. They have been done with a program written by Andrzej Lecki. 
The author is very grateful to him and Zbigniew Szafraniec for giving him this program. 

In this paper, ``smooth" means ``of class at least $C^1$". The ball in $\mathbb{R}^n$ centered at the origin of radius $r$ will be
denoted by $B_r^n$ and $S^r_{n-1}$ is its boundary. If $x $ is in $\mathbb{R}^n$ then $\vert x \vert$ denotes its usual euclidian
norm. Moreover,
we will use the following notations : if $F=(F_1,\ldots,F_k) : \mathbb{R}^n \rightarrow \mathbb{R}^k$, $0< k \le n$, is a smooth mapping then $DF$ is its
Jacobian matrix and $\frac{\partial(F_1,\ldots,F_k)}{\partial(x_{i_1},\ldots,x_{i_k})}$ is the determinant of the
following $k \times k$ minors of $DF$ :
$$\left( \begin{array}{ccc}
{F_1}_{x_{i_1}} & \cdots & {F_1}_{x_{i_k}} \cr
\vdots & \ddots & \vdots \cr
{F_k}_{x_{i_1}} & \cdots & {F_k}_{x_{i_k}} \cr
\end{array} \right).$$

The author is grateful to Jean-Paul Brasselet and David Trotman for their careful reading
of this manuscript and for their remarks and comments. The reader interested in vector fields and $1$-forms
on singular spaces can refer to the monograph [BSS], which gives a detailed account of all the results in this topic.

\section{1-forms on smooth manifolds}
In this section, we recall some well-known facts and results about 1-forms on manifolds. Let $V$ be a smooth manifold of
dimension $n$ and let $\omega$ be a smooth 1-form on $V$. This means that $\omega$ assigns to each point $x$ in $V$ an element
in $(T_x V)^*$, the dual space of $T_x V$. A point $P$ in $V$ is a zero (or a singular point) of $\omega$ if $\omega(P)=0$. We
remark that if $n=0$ then each point in $V$ is a singular point of $\omega$. 

If $P$ is an isolated zero of $\omega$, we can define the index of $\omega$ at $P$. If dim $V=0$, this index is defined to be
$1$. If dim $V >0$, let $\phi : U \subset \mathbb{R}^n \rightarrow V$ be a local parametrization of $V$ at $p$. We can assume
that $\phi(0)=P$. Then the 1-form $\phi^*\omega$ has an isolated zero at $0$. Since ${\mathbb{R}^n}^*$ is isomorphic to
$\mathbb{R}^n$, $\phi^* \omega$ can be viewed as a mapping from $U\subset \mathbb{R}^n$ to $\mathbb{R}^n$. The index of $\omega$
at $P$ is defined to be the degree of the mapping $\frac{\phi^* \omega}{\vert \phi^* \omega\vert} : S_\varepsilon^{n-1} \rightarrow
S^{n-1}$, where $S_\varepsilon^{n-1}$ is a sphere centered at the origin of radius $\varepsilon$ such that $0$ is the only zero of
$\phi^* \omega$ in $B_\varepsilon^n$. Of course, this definition does not depend on the choice of the parametrization. We will
denote by Ind$_{PH}(\omega,P,V)$ this index. When dim $V>0$, we say that $P$ is a non-degenerate zero (or singular point) of $\omega$
if det $D\phi^*\omega \not= 0$. In this case, Ind$_{PH}(\omega,P,V)$ is the sign of the determinant of $D\phi^*\omega(0)$. 
A 1-form $\omega$ on $V$ is non-degenerate if all its zeros are non-degenerate. The set of non-degenerate 1-forms on $V$ is
dense in the set of 1-forms on $V$. If $V$ is compact and $\omega$ is a 1-form on $V$ with a finite number of zeros $P_1,\ldots,P_k$ then the Poincar\'e-Hopf
theorem asserts that $\chi(V)= \sum_{i=1}^k\hbox{Ind}_{PH}(\omega,P_i,V)$. 

If $W$ is a submanifold of $V$ then a 1-form $\omega$ naturally restricts to a 1-form $\omega_{\vert W}$ defined on $W$ in
the following way : for each $x \in W$, $\omega_{\vert W}(x)=\omega(x)_{\vert T_x W}$. We will denote by Ind$_{PH}(\omega,P,W)$
the index Ind$_{PH}(\omega_{\vert W},P,W)$ if $P$ is a zero of $\omega_{\vert W}$. 

From now on, we assume that $\omega$ is a 1-form on an open set $U \subset \mathbb{R}^n$ given by :
$\omega=a_1dx_1+\cdots+a_ndx_n,$
where the $a_i$'s are smooth functions on $U$. Let $V$ be a submanifold of dimension $n-k$ in $U$ and let $P$ be a point in $V$.
We assume that around $P$, $V$ is defined by the vanishing of $k$ smooth functions $f_1,\ldots,f_k$ and that 
$\frac{\partial(f_1,\ldots,f_k)}{\partial(x_1,\ldots,x_k)}(P)\not= 0$. For $j \in \{k+1,\ldots,n\}$, let $m_j$ be defined by :
$$m_j =\left\vert \begin{array}{cccc}
a_1 & \cdots & a_k & a_j \cr
\frac{\partial f_1}{\partial x_1} & \cdots & \frac{\partial f_1}{\partial x_k} & \frac{\partial f_1}{\partial x_j}\cr
\vdots & \ddots & \vdots & \vdots \cr
\frac{\partial f_k}{\partial x_1} & \cdots & \frac{\partial f_k}{\partial x_k} & \frac{\partial f_k}{\partial x_j}\cr
\end{array}
\right\vert.$$
The following lemma tells us when $P$ is a zero of $\omega_{\vert V}$ and, in case it is non-degenerate, gives a way to compute 
Ind$_{PH}(\omega,P,V)$.
\begin{lemma}
The point $P$ is a zero of $\omega_{\vert V}$ if and only if for each $j \in \{k+1,\ldots,n\}$, $m_j(P)=0$. Furthermore it is
non-degenerate if and only if :
$$\frac{\partial(f_1,\ldots,f_k,m_{k+1},\ldots,m_n)}{\partial(x_1,\ldots,x_n)}(P)\not= 0.$$
In this case, 
$$\displaylines{
\quad \hbox{\em Ind}_{PH}(\omega,P,V)= \hfill \cr
\hfill \hbox{\em sign} 
\left( (-1)^{k(n-k)} \frac{\partial(f_1,\ldots,f_k)}{\partial(x_1,\ldots,x_k)}(P)^{n-k+1} 
\frac{\partial(f_1,\ldots,f_k,m_{k+1},\ldots,m_n)}{\partial(x_1,\ldots,x_n)}(P) \right).  \cr
}$$
\end{lemma}
{\it Proof.} The proof is given in [Sz3, p.348-351] in details when $\omega$ is the differential of a function $g$. It
also works in the general case. $\hfill \Box$

Let $(x,\lambda)=(x_1,\ldots,x_n,\lambda_1,\ldots,\lambda_k)$ be a coordinate system in $\mathbb{R}^n \times \mathbb{R}^k$ and
let $H : U \times \mathbb{R}^k \rightarrow \mathbb{R}^n \times \mathbb{R}^k$ be the map given by :
$$\displaylines{
\qquad H(x,\lambda)=\Big( a_1(x)+\sum_{i=1}^k \lambda_i \frac{\partial f_i}{\partial x_1}(x),\ldots,
a_n(x)+\sum_{i=1}^k \lambda_i \frac{\partial f_i}{\partial x_n}(x), \hfill \cr
\hfill  f_1(x),\ldots,f_k(x)  \Big). \qquad \cr
}$$
The following lemma also characterizes a zero of $\omega_{\vert V}$ and computes its index.
\begin{lemma}
The point $P$ is a zero of $\omega_{\vert V}$ if and only if there is a (uniquely determined) point $\lambda \in \mathbb{R}^k$
such that $H(P,\lambda)=0$. Furthermore it is non-degenerate if and only if {\em det}$[DH(P,\lambda)]\not= 0$. In this case, 
$$\hbox{\em Ind}_{PH}(\omega,P,V)=\hbox{\em sign} 
\left( (-1)^k \hbox{\em det}[DH(P,\lambda)] \right).$$
\end{lemma}
{\it Proof.} The lemma is proved carefully when $\omega$ is the differential of a function in [Sz2, Section 1]. The same method
can be applied in the general situation. $\hfill \Box$

\section{A Poincar\'e-Hopf theorem for manifolds with corners}
In this section, we give a version of the Poincar\'e-Hopf theorem for 1-forms defined on a manifold with corners.
First we recall some basic facts about manifolds with corners. Our reference is [Ce].  A manifold
with corners $M$ is defined by an atlas of charts mo\-delled on open subsets of $\mathbb{R}^n_+$.
We write $\partial M$ for its boundary.  
We will make the additional assumption that the boundary is partitioned into pieces $\partial_i M$,
themselves manifolds with corners, such that in each chart, the intersections with the coordinate
hyperplanes $x_j=0$ correspond to distinct pieces $\partial_i M$ of the boundary. For any set $I$ of suffices, we
write $\partial_I M=\cap_{i \in I} \partial_i M$ and we make the convention that 
$\partial_{\emptyset} M = M\setminus \partial M$. 

Any $n$-manifold $M$ with corners can be embedded in a $n$-manifold $M^+$ without boundary so that
the pieces $\partial_i M$ extend to submanifolds $\partial_i M^+$ of codimension 1 in $M^+$.

Let $M$ be a manifold with corners and let $\omega$ be a smooth 1-form on $M^+$.
\begin{definition}
We say that $P$ in $M$ is a zero (or singular point) of $\omega$ on $M$ if it is a zero of a form $\omega_{\vert
\partial_I M^+}$. A zero $P$ of $\omega$ on $M$ is a correct point if, taking $I(P)=\{i \ \vert \ P \in \partial_i M\}$,
$P$ is a zero of $\omega_{\vert \partial_{I(P)} M^+}$ but not a zero of $\omega_{\vert \partial_J M^+}$ for
any proper subset $J$ of $I(P)$. 

A zero $P$ of $\omega$ on $M$ is a non-degenerate correct zero if it is a correct zero of $\omega$
on $M$ and if $P$ is a non-degenerate zero of $\omega_{\vert \partial_{I(P)}M^+}$.
\end{definition}
Note that a $0$-dimensional corner
point $P$ is always a zero because in this case $\partial_{I(P)} M^+=\{P\}$, which is a $0$-dimensional manifold.

\begin{definition}
We say that $\omega$ is a correct (resp. correct non-degenerate) 1-form on $M$ if it admits only correct (resp. correct
non-degenerate) zeros on $M$.
\end{definition}

\begin{proposition}
The set of 1-forms defined on $M^+$ which are correct non-degenerate on $M$ is dense in the set of 1-forms on $M^+$.
\end{proposition}
{\it Proof.} This is clear because there is a finite number of pieces $\partial_I M^+$. $\hfill \Box$

The index Ind$_{PH}(\omega,P,M)$ of $\omega$ on $M$ at a correct zero $P$ is defined to be
Ind$_{PH}(\omega,P,\partial_{I(p)}M^+)$. If $P$ is a correct zero of $\omega$ on $M$, $i \in I(P)$, and $J$ is formed from $I(P)$ by deleting $i$, then
in a chart at $P$ with $\partial_J M^+$ mapping to $\mathbb{R}^p_+$ and $\partial_{I(P)} M$ to the subset $\{x_1=0\}$, the form
$\omega$ on $\partial_J M^+$ has no zeros but its restriction to $\{x_1=0\}$ has one at $P$. Hence $\langle \omega(P),dx_1(P)
\rangle \not=0$, where here the scalar product is considered in ${\mathbb{R}^p}^*$.  

\begin{definition}
We say that $\omega$ is inward at $P$, if for each $i \in I(P)$, we have $\langle \omega(P),dx_1(P) \rangle >0$.
\end{definition}

\begin{remark}
By our convention, if $I(P)= \emptyset$, then $\omega$ is inward at $P$.
\end{remark}

\begin{theorem}
If $M$ is compact and $\omega$ is correct then :
$$\displaylines{
\qquad\chi(M)=\sum \left\{ \hbox{\em Ind}_{PH}(\omega,P,M) \  \vert \  P\ a\  correct\  zero\ of\ \omega \right. \hfill \cr
\hfill \left.   which\ is\ inward\
at\ P  \right\}. \qquad \cr
}$$
\end{theorem}
{\it Proof.} Let us prove it first when $M$ is a manifold with boundary. In this case, it follows from Arnol'd's results
[Ar] mentioned in the introduction. To see this, we just have to relate the index Ind$_{PH}(\omega,M,P)$ when $P$ belongs to
the boundary to the index $i_+(P)$ defined by Arnol'd. We can work in a local chart and assume that $P=0$ in
$\mathbb{R}^n$, that $M=\{ x \in \mathbb{R}^n \ \vert \ x_1 \ge 0\}$ and that $\omega=a_1 dx_1+\cdots+a_n dx_n$.
Then we have (see [Ar,p.7]) :
$$i_+(P)=\frac{1}{2}\left(\hbox{Ind}_{PH}(V,0,\mathbb{R}^n)+\hbox{Ind}_{PH}(V_1,0,\mathbb{R}^n)+
\hbox{Ind}_{PH}(V_0,0,\mathbb{R}^n) \right),$$ where $V$, $V_1$ and $V_0$ are the following
vector fields :
$$V=x_1 a_1 \frac{\partial}{\partial x_1}+a_2 \frac{\partial}{\partial x_2} +\cdots+ a_n
\frac{\partial}{\partial x_n},$$
$$V_1= a_1 \frac{\partial}{\partial x_1}+ a_2 \frac{\partial}{\partial x_2} +\cdots+ a_n
\frac{\partial}{\partial x_n},$$
$$V_0= a_2 \frac{\partial}{\partial x_2} +\cdots+ a_n \frac{\partial}{\partial x_n} \hbox{ on }\{x_1=0\}.$$
Here Ind$_{PH}(V_1,0,\mathbb{R}^n)=0$ since $a_1(P) \not= 0$ and $$\hbox{Ind}_{PH}(V_0,0,\{x_1=0\}) =
\hbox{Ind}_{PH}(\omega,P,M).$$ Furthermore, if $a_1(P) >0$
then Ind$_{PH}(V,0,\mathbb{R}^n)$ is Ind$_{PH}(V_0,0,\{x_1=0\})$ and if $a_1(P)<0$ then it is
$-\hbox{Ind}_{PH}(V_0,0,\{x_1=0\})$. Hence
$i_+(P)=\hbox{Ind}_{PH}(\omega,P,M)$ if $P$ is inward and $i_+(P)=0$ if $P$ is not inward.

Now we suppose that $M$ is a manifold with corners and that $\omega$ is a correct non-degenerate 1-form on $M$. Let us
denote by $Q_1,\ldots,Q_s$ the zeros of $\omega$ lying in $\partial_\emptyset M$ and by $P_1,\ldots,P_r$ those lying in
$\partial M$. Let $h : M \rightarrow \mathbb{R}$ be a carpeting function for $\partial M$ (see the appendix of Douady and
H\'erault in [BoSe]) and let
$\varepsilon'>0$ be a small regular value of $h$ such that $\chi(M)=\chi(M \cap \{h \ge \varepsilon\})$ 
 and $Q_1,\ldots,Q_s$ lie in $M \cap \{ h > \varepsilon\}$, for all
$\varepsilon$ with $0 < \varepsilon \le \varepsilon'$.  Let us study the situation around a point $P_i$. We can find a 
chart $x=(x_1,\ldots,x_n)$ centered at $P_i$ such that in this chart $h$ is the
function $x_1\cdots x_k$ and $\partial_{I(P_i)} M^+$ is the manifold $\{x_1=\cdots=x_k=0\}$ and $M$ is $\{x_1 \ge 0,\ldots,x_k
\ge 0\}$. If we write $\omega=a_1dx_1+\cdots+a_n dx_n$ then $a_{k+1}(P_i)=\cdots=a_n(P_i)=0$ and $a_j(P_i) \not= 0$ for $j \in
\{1,\ldots,k\}$ because $P_i$ is a correct zero of $\omega$. Let $\omega_i$ be the 1-form defined in this chart by :
$$\omega_i(x)=\sum_{j=1}^ka_j(P_i)dx_j+\sum_{j=k+1}^n a_j(x)dx_j.$$
Gluing the initial form $\omega$ with the forms $\omega_i$, we can construct a new form $\tilde{\omega}$ on $M$ with the
following properties :
\begin{itemize}
\item $\tilde{\omega}$ is a correct non-degenerate 1-form on $M$,
\item $\tilde{\omega}=\omega_i$ in a neighborhood of $P_i$,
\item $\tilde{\omega}$ has exactly the same zeros as $\omega$ and the same inward zeros as $\omega$,
\item if $X$ is one of these zeros then Ind$_{PH}(\tilde{\omega},X,M)=\hbox{Ind}_{PH}(\omega,X,M)$.
\end{itemize}
For $\varepsilon >0$ small enough, $\tilde{\omega}$ is clearly a correct 1-form on $\{ h \ge \varepsilon\}$. It is also
non-degenerate for, otherwise we could find a sequence of points $X_k$ such that $h(X_k)=\frac{1}{k}$ and $X_k$ is a degenerate
zero of $\tilde{\omega}_{\vert \{h=\frac{1}{k}\}}$. We can assume that $(X_k)$ tends to a point
$X_0$ in $\{h=0\}$. Using local coordinates around $X_0$, it is easy to see that $X_0$ is a zero of $\tilde{\omega}$, hence
there exists $i \in \{1,\ldots,r\}$ such that $X_0=P_i$. Using Lemma 2.2 and the expression of $\tilde{\omega}$ in a local
chart around $P_i$, we see that $P_i$ is a degenerate zero of $\tilde{\omega}$, which is impossible.  

Let us denote by $P_1,\ldots,P_u$, $u \le r$, the inward critical points of $\tilde{\omega}$. With the expression of $h$
and $\tilde{\omega}$ in local coordinates around $P_i$, it is not difficult to see that each $P_i$, $i \in
\{1,\ldots,u\}$, gives rise to exactly one inward critical point $P_i^\varepsilon$ of $\tilde{\omega}$ on $ \{h \ge
\varepsilon\}$. Furthermore, using Lemma 2.2 and making some computations of determinants, we find that this critical point
$P_i^\varepsilon$ is non-degenerate and has the same index as $\tilde{\omega}$ at $P_i$. Applying the Poincar\'e-Hopf
theorem for manifolds with boundary, we get the result for a correct non-degenerate 1-form. If the form is correct but
admits degenerate zeros, we perturb it around its degenerate zeros and apply the previous case. $\hfill \Box$

\begin{remark}
Since a manifold with corners is a Whitney stratified set, it would be interesting to deduce the above result from
Poincar\'e-Hopf theorems for stratified sets like Theorem 1 in [EG5], Theorem 5.4 in [KT], Theorem 6.2.2
in [Sc3] or Theorem 2 in [Si].
\end{remark}

\section{The radial index of a 1-form}
The notion of radial index was defined by Ebeling and Gusein-Zade for 1-forms on real analytic sets in [EG5]. This notion is inspired
by the work of M.H Schwartz on radial vector fields on singular analytic varieties.
Here we recall the definition of the radial index of a $1$-form but in the more general setting of closed subanalytic sets.

Let $X \subset \mathbb{R}^n$ be a closed subanalytic set equipped with a Whitney stratification $\{S_\alpha\}_{\alpha \in
\Lambda}$. Let $\omega$ be a continuous 1-form defined on $\mathbb{R}^n$. We say that a point $P$ in $X$ is a zero (or a
singular point) of $\omega$
on $X$ if it is a zero of $\omega_{\vert S}$, where $S$ is the stratum that contains $P$. In the sequel, we will define the radial
index of $\omega$ at $P$, when $P$ is an isolated zero of $\omega$ on $X$. We can assume that $P=0$ and we denote by $S_0$ the stratum that contains $0$.

\begin{definition}
A 1-form $\omega$ is radial on $X$ at $0$ if, for an arbitrary non-trivial subanalytic arc $\varphi : [0,\nu[
\rightarrow X$ of class $C^1$, the value of the form $\omega$ on the tangent vector $\dot{\varphi}(t)$ is positive for $t$
small enough.
\end{definition}

Let $\varepsilon >0$ be small enough so that in the closed ball $B_\varepsilon^n$ of radius $\varepsilon$ centered at $0$ in
$\mathbb{R}^n$, the 1-form has no singular points on $X \setminus \{0\}$. Let $V_0,\ldots,V_q$ be the strata that contain $0$
in their closure. Following Ebeling and Gusein-Zade, there exists a 1-form $\tilde{\omega}$ on
$\mathbb{R}^n$ such that : 
\begin{enumerate}
\item The 1-form $\tilde{\omega}$ coincides with the 1-form $\omega$ on a neighborhood of $S_\varepsilon^{n-1}=\partial
B_\varepsilon^n$.
\item The 1-form $\tilde{\omega}$ is radial on $X$ at the origin.
\item In a neighborhood of each zero $Q \in X \cap B_\varepsilon^n \setminus \{0\}$, $Q \in V_i$, dim$ V_i=k$, the
1-form $\tilde{\omega}$ looks as follows. There exists a local subanalytic diffeomorphism $h :(\mathbb{R}^n,\mathbb{R}^k,0)
\rightarrow (\mathbb{R}^n,V_i,Q)$ such that $h^* \tilde{\omega}=\pi_1^* \tilde{\omega}_1+\pi_2^*\tilde{\omega}_2$ where
$\pi_1$ and $\pi_2$ are the natural projections $\pi_1 : \mathbb{R}^n \rightarrow \mathbb{R}^k$ and $\pi_2 : \mathbb{R}^n
\rightarrow \mathbb{R}^{n-k}$, $\tilde{\omega}_1$ is a 1-form on a neighborhood of $0$ in $\mathbb{R}^k$ with an isolated 
zero at the origin and $\tilde{\omega}_2$ is a radial 1-form on $\mathbb{R}^{n-k}$ at $0$.
\end{enumerate}

\begin{definition}
The radial index $\hbox{\em Ind}_{Rad}(\omega,0,X)$ of the 1-form $\omega$ on $X$ at $0$ is the sum :
$$1+\sum_{i=1}^q \sum_{Q \vert  \tilde{\omega}_{\vert V_i}(Q)=0} \hbox{\em Ind}_{PH}(\tilde{\omega},Q,V_i),$$
where the sum is taken over all zeros of the 1-form $\tilde{\omega}$ on $(X\setminus \{0\}) \cap B_\varepsilon$. 
If $0$ is not a zero of $\omega$ on $X$, we put $\hbox{\em Ind}_{Rad}(\omega,0,X)=0$.
\end{definition}
A straightforward corollary of this definition is that the radial index satisfies the law of conservation of number (see 
Remark 9.4.6 in [BSS] or the remark before Proposition 1 in [EG5]).

As in the case of an analytic set, this notion is well defined, i.e it does not depend on the different choices made to
define it. Furthermore, the Poincar\'e-Hopf theorem proved in [EG5] also holds for compact subanalytic sets, with the same
proof. 

\section{The radial index on a manifold with corners}
In this section, we relate the radial index of a 1-form on a manifold with corners to usual
Poincar\'e-Hopf indices of 1-forms.

We work in $\mathbb{R}^n$ with coordinates $(x_1,\ldots,x_n)$. For $1\le k \le n$ and for every
$\epsilon=(\epsilon_1,\ldots,\epsilon_k) \in \{0,1\}^k$, let $\mathbb{R}^n(\epsilon)$ be the following manifold with corners
: $$\mathbb{R}^n(\epsilon)=\left\{ (x_1,\ldots,x_n)\in \mathbb{R}^n \ \vert \ (-1)^{\epsilon_1} x_1 \ge
0,\ldots,(-1)^{\epsilon_k} x_k \ge 0\right\}.$$
Now we consider a smooth 1-form $\Omega=a_1dx_1+\cdots+a_ndx_n$ on $\mathbb{R}^n$. We will denote by $A$ the set
$\{(0,1),(1,0),(1,1)\}$. For every $k \in
\{1,\ldots,n\}$, for every $\underline{\alpha}=((\alpha_1,\beta_1),\ldots,(\alpha_k,\beta_k)) \in A^k$, we define the vector field
$V(\underline{\alpha})$ in the following way :
$$V(\underline{\alpha})=x_1^{\alpha_1}a_1^{\beta_1}\frac{\partial}{\partial x_1}+\cdots+
x_k^{\alpha_k}a_k^{\beta_k}\frac{\partial}{\partial x_k}+a_{k+1}\frac{\partial}{\partial x_{k+1}}+\cdots+
a_n\frac{\partial}{\partial x_n}.$$
We will denote by $\underline{1}$ the element $((1,1),\ldots,(1,1))$.
\begin{proposition}
The form $\Omega$ has an isolated zero at $0$ on $\mathbb{R}^n(\epsilon)$ for every $\epsilon \in \{0,1\}^k$ if and
only if the vector field $V(\underline{1})$ has an isolated zero at the origin.
\end{proposition}
{\it Proof.} The form $\Omega$ has an isolated zero at $0$ on $\mathbb{R}^n(\epsilon)$ for every $\epsilon \in
\{0,1\}^k$ if and only if for every $\underline{\alpha} \in A^k$, the vector field
$V(\underline{\alpha})$ has an isolated zero at the origin. This is equivalent to the fact that $V(\underline{1})$ has an
isolated zero. $\hfill \Box$

From now on, we assume that $V(\underline{1})$ has an isolated zero at the origin. Since $\mathbb{R}^n(\epsilon)$ is
clearly a subanalytic set and $\Omega$ has an isolated zero at $0$ on $\mathbb{R}^n(\epsilon)$, the radial index of
$\Omega$ on $\mathbb{R}^n(\epsilon)$ at the origin is well-defined. 
For each $r>0$, $B_r^n(\epsilon)=B_r^n \cap \mathbb{R}^n(\epsilon)$ and $S_r^{n-1}(\epsilon)=S_r^{n-1} \cap
\mathbb{R}^n(\epsilon)$ are manifolds with corners. Let $\tilde{\Omega}_r$ be a small perturbation of $\Omega$ such that
$\tilde{\Omega}_r$ is correct on $B_r^{n}(\epsilon)$. This implies that $\tilde{\Omega}_r$ is also correct on 
$S_r^{n-1}(\epsilon)$. 
In this situation, we can relate Ind$_{Rad}(\Omega,0,\mathbb{R}^n(\epsilon))$ to the zeros of $\tilde{\Omega}_r$ on
$S_r^{n-1}(\epsilon)$.
\begin{lemma}
Let $\{P_i\}$ be the set of inward zeros of $\tilde{\Omega}_r$
on $B_r^n(\epsilon)$ lying in $S_r^{n-1}$. We have :
$$\hbox{\em Ind}_{Rad}(\Omega,0,\mathbb{R}^n(\epsilon))=1- \sum_i
\hbox{\em Ind}_{PH}(\tilde{\Omega}_r,P_i,S_r^{n-1}(\epsilon)).
$$
\end{lemma}
{\it Proof.} Let us consider first the case when $0$ is a zero of $\Omega$ on $\mathbb{R}^n(\epsilon)$. As a manifold with corners, the set $\mathbb{R}^n(\epsilon)$ has a natural Whitney stratification. Hence we can write $\mathbb{R}^n(\epsilon)=
\cup_{i=0}^q V_i$, where $0 \in V_0$. Let $\tilde{\omega}$ be a 1-form on
$\mathbb{R}^n$ such that :
\begin{enumerate}
\item the 1-form $\tilde{\omega}$ coincides with the 1-form $\Omega$ on a neighborhood of $S_r^{n-1}$,
\item the 1-form $\tilde{\omega}$ is radial in $\mathbb{R}^n(\epsilon)$ at the origin,
\item in a neighborhood of each zero $Q \in \mathbb{R}^n(\epsilon) \cap B_r \setminus \{0\}$, $Q \in V_i$, dim$ V_i=k$, the
1-form $\tilde{\omega}$ looks as follows. There exists a local diffeomorphism $h :(\mathbb{R}^n,\mathbb{R}^k,0)
\rightarrow (\mathbb{R}^n,V_i,Q)$ such that $h^* \tilde{\omega}=\pi_1^* \tilde{\omega}_1+\pi_2^*\tilde{\omega}_2$ where
$\pi_1$ and $\pi_2$ are the natural projections $\pi_1 : \mathbb{R}^n \rightarrow \mathbb{R}^k$ and $\pi_2 : \mathbb{R}^n
\rightarrow \mathbb{R}^{n-k}$, $\tilde{\omega}_1$ is the germ of a 1-form on $(\mathbb{R}^k,0)$ with an isolated zero at
the origin and $\tilde{\omega}_2$ is a radial 1-form on $(\mathbb{R}^{n-k},0)$.
\end{enumerate}
We have :
$$\hbox{Ind}_{Rad}(\Omega,0,\mathbb{R}^n(\epsilon))=1+\sum_{i=1}^q \sum_{Q \vert  \tilde{\omega}_{\vert V_i}(Q)=0} 
\hbox{Ind}_{PH}(\tilde{\omega},Q,V_i).$$
Let us focus on the situation around a zero $Q$ of $\tilde{\omega}$ on $\mathbb{R}^n(\epsilon)$. It is not a correct zero in
the sense of Section 3, because the form $\tilde{\omega}_2$ that appears in the point (3) above is radial. However, if we replace
$\tilde{\omega}_2$ by a small perturbation $\tilde{\omega}_2'=\tilde{\omega}_2 -u_1dx_1-\cdots-u_{n-k}dx_{n-k}$ where
$u_i\not=0$ for each $i \in \{1,\ldots,n-k\}$, then the 1-form $\tilde{\omega}'={h^{-1}}^*(\pi_1^* \tilde{\omega}_1 +
\pi_2^*\tilde{\omega}_2')$ is a correct 1-form in the neighborhood of $Q$ in $\mathbb{R}^n(\epsilon)$. Furthermore it
admits exactly one inward correct singular point $\tilde{Q}$ in the neighborhood of $Q$ which lies in a stratum $V_j$ such
that dim $V_j \ge \hbox{ dim} V_i$ and Ind$_{PH} (\tilde{\omega}',\tilde{Q},V_j)$ is equal to 
Ind$_{PH} (\tilde{\omega},Q,V_i)$. Let $r'$, $0<r'<r$ be such that the points $Q$'s above lie in $\{r' <\vert x \vert <r\}$. 
We can construct a 1-form $\tilde{\omega}'$ on $\mathbb{R}^n$ close to $\tilde{\omega}$ 
such that :
\begin{enumerate}
\item $\tilde{\omega}'$ is a correct 1-form on $\mathbb{R}^n(\epsilon) \cap \{ r' \le \vert x \vert \le r\}$,
\item $\tilde{\omega}'$ coincides with $\tilde{\Omega}_r$ in a neighborhood of $S_r^{n-1}$,
\item 
$$
 \sum_{i=1}^q \sum_{Q \vert  \tilde{\omega}_{\vert V_i}(Q)=0} \hbox{Ind}_{PH}(\tilde{\omega},Q,V_i)= 
\sum_j \hbox{Ind}_{PH} (\tilde{\omega}',Q'_j, \mathbb{R}^n(\epsilon)),$$
where $\{Q'_j\}$ is the set of inward correct zeros of $\tilde{\omega}'$ on $\mathbb{R}^n(\epsilon)$ 
in $\{r' < \vert x \vert < r \}$. 
\item the zeros of $\tilde{\omega}'$ lying in $S_{r'}^{n-1}$ are inward for $\mathbb{R}^n(\epsilon) \cap \{ r' \le \vert x \vert
\le r\}$.
\end{enumerate}

If we denote by $\{S_l\}$ the set of inward correct zeros of $\tilde{\omega}'$ on $\mathbb{R}^n(\epsilon) \cap \{ r' \le \vert x \vert
\le r\}$ such that 
$\vert S_l \vert=r'$ then, by the Poincar\'e-Hopf theorem (Theorem 3.6), we get :
$$1=\chi(\mathbb{R}^n(\epsilon) \cap \{r' \le \vert x \vert \le r\})= 
\sum_l \hbox{Ind}_{PH}(\tilde{\omega}',S_l,\mathbb{R}^n(\epsilon) \cap \{r' \le \vert x \vert \le r\}) $$
$$-1+\hbox{Ind}_{Rad}(\Omega,0,\mathbb{R}^n(\epsilon)) + 
\sum_i \hbox{Ind}_{PH}(\tilde{\Omega}_r,P_i,\mathbb{R}^n(\epsilon) \cap \{r' \le \vert x \vert \le r\}). $$
Since $\tilde{\omega}'$ is correct on $\mathbb{R}^n(\epsilon) \cap \{ r' \le \vert x \vert \le r\}$, it is also correct on
$\mathbb{R}^n(\epsilon) \cap S_{r'}^{n-1}$. Applying the Poincar\'e-Hopf theorem and using point (4) above, we obtain :
$$\sum_l\hbox{Ind}_{PH}(\tilde{\omega}',S_l,\mathbb{R}^n(\epsilon) \cap \{r' \le \vert x \vert \le r\}) =$$
$$\sum_l \hbox{Ind}_{PH}(\tilde{\omega}',S_l,S_{r'}(\epsilon)) =\chi(S_{r'}(\epsilon))=1. $$
It is easy to conclude because for each $i$, we have :
$$\hbox{Ind}_{PH}(\tilde{\Omega}_r,P_i,\mathbb{R}^n(\epsilon) \cap \{r' \le \vert x \vert \le r\})=\hbox{Ind}_{PH}
(\tilde{\Omega}_r,P_i,S_r^{n-1}(\epsilon)).$$ 
When $0$ is not a zero of $\Omega$ on $\mathbb{R}^n(\epsilon)$, we can write : 
$$ 1=\chi(B_r^n(\epsilon))= \sum_i
\hbox{Ind}_{PH}(\tilde{\Omega}_r,P_i,S_r^{n-1}(\epsilon)).$$
The result is proved because Ind$_{Rad}(\Omega,0,\mathbb{R}^n(\epsilon))=0$. 
$\hfill \Box$

Note that this characterization of the radial index is very similar to the definition of the index at an
isolated zero or virtual zero  of a vector field on a radial manifold complex of  King and Trotman ([KT], Definition 5.5).
Now let $\tilde{\Omega'}_r$ be a small perturbation of $\Omega$ such that $\tilde{\Omega'}_r$
is correct on $B_r^n(\epsilon)$. We can relate Ind$_{Rad}(\Omega,0,\mathbb{R}^n(\epsilon))$ to the zeros of
$\tilde{\Omega'}_r$ on $B_r^n(\epsilon)$.
\begin{lemma}
Let $\{Q_j\}$ be the set of inward zeros of 
$\tilde{\Omega'}_r$ on $B_r^n(\epsilon)$ lying in $\{ \vert x \vert < r \}$. We have :
$$\hbox{\em Ind}_{Rad}(\Omega,0,\mathbb{R}^n(\epsilon)) =
\sum_j \hbox{\em Ind}_{PH}(\tilde{\Omega'}_r,Q_j,B_r^n(\epsilon)).$$
\end{lemma}
{\it Proof.} If $\{R_l\}$ is the set of inward zeros
 of $\tilde{\Omega'}_r$ on $B_r^n(\epsilon)$ then, by the Poincar\'e-Hopf theorem, we have :
$$1=\chi(B_r^n(\epsilon)) =\sum_l \hbox{Ind}_{PH}(\tilde{\Omega'}_r,R_l,B_r^n(\epsilon)).$$ 
Now we can decompose $\{R_l\}$ into $\{R_l\}=\{Q_j\} \sqcup \{P_i\}$ where the $P_i$'s are the 
inward zeros of $\tilde{\Omega'}_r$ on $B_r^n(\epsilon)$ lying in $S_r^{n-1}$.
By the previous lemma,
$$\hbox{Ind}_{Rad}(\Omega,0,\mathbb{R}^n(\epsilon))=1-\sum_i \hbox{Ind}_{PH} (\tilde{\Omega'}_r,P_i,S_r^{n-1}(\epsilon)).$$
Summing these two equalities gives the result. $\hfill \Box$

We can state the main result of this section.
\begin{theorem}
Assume that $V(\underline{1})$ has an isolated zero at the origin. For every $\epsilon=(\epsilon_1,\ldots,\epsilon_k) \in
\{0,1\}^k$, we have :
$$\hbox{\em Ind}_{Rad}(\Omega,0,\mathbb{R}^n(\epsilon))=\frac{1}{2^k} (-1)^{\vert \epsilon \vert}\sum_{\underline{\alpha} \in A^k}
(-1)^{[\epsilon \cdot \underline{\alpha}]} \hbox{\em Ind}_{PH}(
V(\underline{\alpha}),0,\mathbb{R}^n),$$
where $\vert \epsilon \vert =\sum_{i=1}^k \epsilon_i$ and $[\epsilon \cdot \underline{\alpha}]=
\sum_{i=1}^k \epsilon_i(\alpha_i+\beta_i)$.
\end{theorem}
{\it Proof.} We will prove this theorem by induction on $k$. Let us assume first that $k=1$. Let
$\tilde{\Omega}=\tilde{a}_1dx_1+\cdots+\tilde{a}_ndx_n$ be a small perturbation of $\Omega$ such that
$\tilde{\Omega}$ is correct and non-degenerate on $B_r^n(0)$ and $B_r^n(1)$ for $r$ small. Let
$\tilde{V}((0,1))$, $\tilde{V}((1,0))$ and $\tilde{V}((1,1))$ be the following vector fields :
$$\tilde{V}((0,1))=\tilde{a}_1\frac{\partial}{\partial x_1}+\tilde{a}_2\frac{\partial}{\partial x_2}+\cdots+
\tilde{a}_n\frac{\partial}{\partial x_n},$$
$$\tilde{V}((1,0))=x_1\frac{\partial}{\partial x_1}+\tilde{a}_2\frac{\partial}{\partial x_2}+\cdots+
\tilde{a}_n\frac{\partial}{\partial x_n},$$
$$\tilde{V}((1,1))=x_1\tilde{a}_1\frac{\partial}{\partial x_1}+\tilde{a}_2\frac{\partial}{\partial x_2}+\cdots+
\tilde{a}_n\frac{\partial}{\partial x_n}.$$
For $r$ small enough, for $\underline{\alpha} \in \{(0,1),(1,0),(1,1)\}$, the degree of the mapping
$\frac{\tilde{V}(\underline{\alpha})}{\vert \tilde{V}(\underline{\alpha}) \vert} : S_r^{n-1} \rightarrow S^{n-1}$ is equal to
Ind$_{PH} (V(\underline{\alpha}),0,\mathbb{R}^n)$. Furthermore, the zeros of $\tilde{V}(\underline{\alpha})$ inside $B_r^n$ are all
non-degenerate by our assumption on $\tilde{\Omega}$. Using this characterization of Ind$_{PH}(V(\underline{\alpha}),0,\mathbb{R}^n)$ and the way to compute
Ind$_{Rad}(\Omega,0,\mathbb{R}^n(0))$ and Ind$_{Rad}(\Omega,0,\mathbb{R}^n(1))$ given in the previous lemma, we find :
$$\displaylines{
\qquad \hbox{Ind}_{Rad}(\Omega,0,\mathbb{R}^n(0)) +\hbox{Ind}_{Rad}(\Omega,0,\mathbb{R}^n(1)) = \hfill \cr
\hfill \hbox{Ind}_{PH}(V((1,0)),0,\mathbb{R}^n) + \hbox{Ind}_{PH}(V((0,1)),0,\mathbb{R}^n), \qquad}
$$
$$\hbox{Ind}_{Rad}(\Omega,0,\mathbb{R}^n(0)) -\hbox{Ind}_{Rad}(\Omega,0,\mathbb{R}^n(1))=
\hbox{Ind}_{PH}(V((1,1)),0,\mathbb{R}^n).$$
This gives the result for $k = 1$. Now assume that $k>1$. Let $\tilde{\Omega}=\tilde{a}_1dx_1+\cdots+\tilde{a}_ndx_n$ be a small perturbation of $\Omega$ such that
$\tilde{\Omega}$ is correct and non-degenerate on $B_r^n(\epsilon)$ for $r$ small enough and for every $\epsilon \in
\{0,1\}^k$. For $\underline{\alpha} \in \{(1,0),(0,1),(1,1)\}^k$, let $\tilde{V}(\underline{\alpha})$ be the vector field
defined by :
$$\tilde{V}(\underline{\alpha})=x_1^{\alpha_1}\tilde{a}_1^{\beta_1}\frac{\partial}{\partial x_1}+\cdots+
x_k^{\alpha_k}\tilde{a}_k^{\beta_k}\frac{\partial}{\partial x_k} +\tilde{a}_{k+1}\frac{\partial}{\partial
x_{k+1}}+\cdots+\tilde{a}_n 
\frac{\partial}{\partial x_n}.$$
As above, if $r$ is small enough then, $\tilde{V}(\underline{\alpha})$ admits only non-degenerate zeros in $B_r^n$ and the
degree of the map $\frac{\tilde{V}(\underline{\alpha})}{\vert \tilde{V}(\underline{\alpha}) \vert} : S_r^{n-1}
\rightarrow S^{n-1}$ is Ind$_{PH} (V(\underline{\alpha}),0,\mathbb{R}^n)$. Let us fix $\epsilon' \in \{0,1\}^{k-1}$ and let
$\epsilon^0=(\epsilon',0)$ and $\epsilon^1=(\epsilon',1)$. Since $\tilde{\Omega}$ is correct and non-degenerate on $B_r^n(\epsilon^0)$ and 
$B_r^n(\epsilon^1)$, it is also correct and non-degenerate on $B_r^n(\epsilon')$
and $B_r^n(\epsilon') \cap \{x_k=0\}$. Counting carefully the
zeros of these vector fields and using the previous lemma, we obtain that :
$$\displaylines{
\qquad \hbox{Ind}_{Rad}(\Omega,0,\mathbb{R}^n(\epsilon^0))+ \hbox{Ind}_{Rad}(\Omega,0,\mathbb{R}^n(\epsilon^1))= \hfill \cr
\hfill \hbox{Ind}_{Rad}(\Omega,0,\mathbb{R}^n(\epsilon'))+\hbox{Ind}_{Rad}(\Omega,0,\mathbb{R}^n(\epsilon')\cap \{x_k=0\}).
\qquad \cr
}$$
Let $\Gamma$ be the 1-form defined by :
$$\Gamma=a_1dx_1+\cdots+a_{k-1}dx_{k-1}+x_ka_kdx_k+a_{k+1}dx_{k+1}+\cdots+a_ndx_n.$$
With the same arguments, we find :
$$\hbox{Ind}_{Rad}(\Omega,0,\mathbb{R}^n(\epsilon^0))- \hbox{Ind}_{Rad}(\Omega,0,\mathbb{R}^n(\epsilon^1))=
\hbox{Ind}_{Rad}(\Gamma,0,\mathbb{R}^n(\epsilon')).$$
It is enough to use the inductive hypothesis to conclude. $\hfill \Box$

We can apply Theorem 5.4 to the differential of an analytic function-germ and use Theorem 2 in [EG5].
\begin{corollary}
Let $f:(\mathbb{R}^n,0) \rightarrow (\mathbb{R},0)$ be an analytic function-germ with an isolated critical point at the
origin. Let $k \in \{1,\ldots,n\}$ and assume that the vector field $\nabla f (\underline{1})$ has an isolated zero at the
origin where $\nabla f$ is the gradient vector field of $f$. Then for every
$\underline{\alpha} \in A^k$, $\nabla f(\underline{\alpha})$ has an isolated zero at the origin and
for $\delta$ such that $0 < \vert \delta \vert \ll r \ll 1$, we have :
$$\displaylines{
\qquad \chi(f^{-1}(\delta) \cap B_r^n \cap \mathbb{R}^n(\epsilon))= \hfill \cr
1-\frac{1}{2^k} (-1)^{\vert \epsilon \vert} \left[
\hbox{\em sign}(-\delta)^{n-k}\sum_{\underline{\alpha} \in A^k \ \vert \ \vert \underline{\alpha} \vert_2  \ even}
(-1)^{[\epsilon \cdot \underline{\alpha}]} \hbox{\em Ind}_{PH}(\nabla f(\underline{\alpha}),0,\mathbb{R}^n) + \right. \cr
\left. \hbox{\em sign}(-\delta)^{n-k+1}\sum_{\underline{\alpha} \in A^k \ \vert \ \vert \underline{\alpha} \vert_2 \ odd}
(-1)^{[\epsilon \cdot \underline{\alpha}]} \hbox{\em Ind}_{PH} (\nabla f(\underline{\alpha}),0,\mathbb{R}^n) \right]. \cr
}$$
where, if $\underline{\alpha}=((\alpha_1,\beta_1),\ldots,(\alpha_k,\beta_k))$ then $\vert \underline{\alpha} \vert_2 =\sum_{i=1}^k
\beta_i$.
\end{corollary}

\begin{remark}
In [Du2], we explained in Section 6 how the Euler-Poincar\'e characteristic of $f^{-1}(\delta) \cap B_r \cap
\mathbb{R}^n(\epsilon))$ can be related to the indices of the vector fields $\nabla
f(\underline{\alpha})$ but we did not give any explicit formula.
\end{remark}

\noindent{\bf Examples}

$ \bullet $ Let $\Omega(x_1,x_2)=(x_1-x_2)dx_1+(x_2^2+x_1x_2)dx_2$.

For $\underline{\alpha}=((1,0),(1,0))$, it is clear that Ind$_{PH}(V(\underline{\alpha}),0,\mathbb{R}^n)=1$. 

For $\underline{\alpha}=((0,1),(1,0))$, we find that Ind$_{PH}(V(\underline{\alpha}),0,\mathbb{R}^n)=1$. 

\noindent Using the program written by Lecki, we can compute the indices of the other $V(\underline{\alpha})$'s. 

For $\underline{\alpha}=((1,0),(0,1))$, Ind$_{PH}(V(\underline{\alpha}),0,\mathbb{R}^n)=0$. 

For $\underline{\alpha}=((0,1),(0,1))$, Ind$_{PH}(V(\underline{\alpha}),0,\mathbb{R}^n)=0$. 

For $\underline{\alpha}=((1,1),(1,0))$, Ind$_{PH}(V(\underline{\alpha}),0,\mathbb{R}^n)=0$. 

For $\underline{\alpha}=((1,0),(1,1))$, Ind$_{PH}(V(\underline{\alpha}),0,\mathbb{R}^n)=1$.

For $\underline{\alpha}=((1,1),(0,1))$, Ind$_{PH}(V(\underline{\alpha}),0,\mathbb{R}^n)=0$. 

For $\underline{\alpha}=((0,1),(1,1))$, Ind$_{PH}(V(\underline{\alpha}),0,\mathbb{R}^n)=1$. 

For $\underline{\alpha}=((1,1),(1,1))$, Ind$_{PH}(V(\underline{\alpha}),0,\mathbb{R}^n)=0$. 

\noindent Applying Theorem 5.4, we obtain :
$$\hbox{Ind}_{Rad}(\Omega,0,\mathbb{R}^2((0,0)))=1, 
\hbox{Ind}_{Rad}(\Omega,0,\mathbb{R}^2((0,1)))=0,$$
$$ \hbox{Ind}_{Rad}(\Omega,0,\mathbb{R}^2((1,0)))=1 ,
\hbox{Ind}_{Rad}(\Omega,0,\mathbb{R}^2((1,1)))=0.$$

$ \bullet $ Let $\Omega(x_1,x_2)=(x_1^2+x_1x_2)dx_1-(2x_1x_2+x_2^2)dx_2$.

For $\underline{\alpha}=((1,0),(1,0))$, it is clear that Ind$_{PH}(V(\underline{\alpha}),0,\mathbb{R}^n)=1$.

\noindent Using the program written by Lecki, we can compute the indices of the other $V(\underline{\alpha})$'s. 

For $\underline{\alpha}=((0,1),(1,0))$, Ind$_{PH}(V(\underline{\alpha}),0,\mathbb{R}^n)=0$. 

For $\underline{\alpha}=((1,0),(0,1))$, Ind$_{PH}(V(\underline{\alpha}),0,\mathbb{R}^n)=0$. 

For $\underline{\alpha}=((0,1),(0,1))$, Ind$_{PH}(V(\underline{\alpha}),0,\mathbb{R}^n)=-2$. 

For $\underline{\alpha}=((1,1),(1,0))$, Ind$_{PH}(V(\underline{\alpha}),0,\mathbb{R}^n)=1$. 

For $\underline{\alpha}=((1,0),(1,1))$, Ind$_{PH}(V(\underline{\alpha}),0,\mathbb{R}^n)=-1$.

For $\underline{\alpha}=((1,1),(0,1))$, Ind$_{PH}(V(\underline{\alpha}),0,\mathbb{R}^n)=0$. 

For $\underline{\alpha}=((0,1),(1,1))$, Ind$_{PH}(V(\underline{\alpha}),0,\mathbb{R}^n)=0$. 

For $\underline{\alpha}=((1,1),(1,1))$, Ind$_{PH}(V(\underline{\alpha}),0,\mathbb{R}^n)=1$. 

\noindent Applying Theorem 5.4, we obtain :
$$\hbox{Ind}_{Rad}(\Omega,0,\mathbb{R}^2((0,0)))=0, 
\hbox{Ind}_{Rad}(\Omega,0,\mathbb{R}^2((0,1)))=0,$$
$$ \hbox{Ind}_{Rad}(\Omega,0,\mathbb{R}^2((1,0)))=0 ,
\hbox{Ind}_{Rad}(\Omega,0,\mathbb{R}^2((1,1)))=0.$$

\section{Condition $(P')$ and its consequences}
The results obtained in this section will be used in the
study of 1-forms on some hypersurfaces with isolated singularities that we will do in the next section.

Let $V=a_1 \frac{\partial}{\partial x_1}+\cdots+a_n\frac{\partial}{\partial x_n}$ be a smooth vector field defined
in a neighborhood of the origin such that $0$ is an isolated zero of $V$. We suppose that $V$ satisfies the following condition
$(P')$ : there exist smooth vector fields $V_2,\ldots,V_n$ defined in the neighborhood of $0$ such that
$V_2(x),\ldots,V_n(x)$ span $V(x)^\perp$ whenever $V(x) \not= 0$ and such that $(V(x),V_2(x),\ldots,V_n(x))$ is a direct basis
of $\mathbb{R}^n$. When $V$ is the gradient vector of a function, Condition $(P')$ coincides with Condition $(P)$
introduced by Fukui and Khovanskii [FK]. 

The following proposition gives necessary and sufficient conditions for the existence of $V_2,\ldots,V_n$.

\begin{proposition}
Let $V$ be a smooth vector field defined in the neighborhood of the origin with an isolated zero at the origin. The following
conditions are equivalent :
\begin{itemize}
\item $V$ satifies Condition $(P')$,
\item one of the following conditions holds :
\begin{itemize}
\item $n=2,4$ or $8$,
\item $n$ is even, $n\not= 2,4,8$, and $\hbox{\em Ind}_{PH}(V,0,\mathbb{R}^n) $ is even,
\item $n$ is odd and $\hbox{\em Ind}_{PH}(V,0,\mathbb{R}^n)=0$.
\end{itemize}
\end{itemize}
\end{proposition}
{\it Proof.} The proof for a gradient vector field is given [FK], Section 1.1. It can be mimicked in the general case.
$\hfill \Box$

Furthermore, when $n=2,4,8$ or $a_1 \ge 0$, it is possible to construct explicitely the vector fields $V_i$ in terms of the
components $a_1,\ldots,a_n$ of $V$ and if $V$ is analytic (resp. polynomial), so are the $V_i$'s. This is explained in [FK],
Section 1.2 for a gradient vector field and works exactly in the same way in the general case.

From now on, we work with a vector field $V=a_1\frac{\partial}{\partial x_1}+\cdots+a_n\frac{\partial}{\partial x_n}$ with an
isolated singularity at the origin, that satisfies Condition $(P')$. Let $X \in S^{n-1}$ and let $W_X$ be the vector field
given by :
$$W_X(x)=\langle V(x),X \rangle \frac{\partial}{\partial x_1}+\langle V_2(x),X \rangle \frac{\partial}{\partial x_2}
+ \cdots +\langle V_n(x),X \rangle  \frac{\partial}{\partial x_n}.$$
\begin{lemma}
The vector field $W_X$ has an isolated zero at the origin.
\end{lemma}
{\it Proof.} If $x\not= 0$ then $(V(x),V_2(x),\ldots,V_n(x))$ is a basis of $\mathbb{R}^n$, so $W_X(x) \not= 0$ because $X
\not= 0$. $\hfill \Box$

\begin{lemma}
For every $X \in S^{n-1}$, $\hbox{\em Ind}_{PH}(W_X,0,\mathbb{R}^n)=\hbox{\em Ind}_{PH}(W_{e_1},0,\mathbb{R}^n)$  where
$e_1=(1,0,\ldots,0)$.
\end{lemma}
{\it Proof.} Let us fix $X \in S^{n-1}$. There exists $A \in SO(n)$ such that $A.X=e_1$. Since $SO(n)$ is arc-connected, $W_X$
and $W_{e_1}$ are homotopic. Furthermore, thanks to Condition $(P')$, we can choose $r$ small enough such that all the
$W_Y$'s, with $Y \in S^{n-1}$, have no zero in
$B_r^n\setminus \{0\}$. Hence the mappings $\frac{W_X}{\vert W_X \vert} : S_r^{n-1} \rightarrow S^{n-1}$ and 
$\frac{W_{e_1}}{\vert W_{e_1} \vert} : S_r^{n-1} \rightarrow S^{n-1}$ are homotopic as well. $\hfill \Box$

Our first aim is to compare $\hbox{Ind}_{PH}(V,0,\mathbb{R}^n)$ and $\hbox{Ind}_{PH}(W_{e_1},0,\mathbb{R}^n)$.

\begin{lemma}
We have :
$$\frac{V}{\vert V \vert}(x)=e_1 \Leftrightarrow \frac{W_{e_1}}{\vert W_{e_1} \vert}(x)=e_1.$$
\end{lemma}
{\it Proof.} If $\frac{V}{\vert V \vert}(x)=e_1$ then $a_1(x)>0$ and $a_i(x)=0$ for $i \in \{2,\ldots,n\}$. Since $a_1(x)>0$,
the family $(V_2'(x),\ldots,V_n'(x))$ is a basis of $V(x)^\perp$ where $V_i'$ is defined by :
$$V_i'=-a_i \frac{\partial}{\partial x_1}+a_1 \frac{\partial}{\partial x_i}.$$
Furthermore $(V(x),V_2'(x),\ldots,V_n'(x))$ is direct. There exists a direct $(n-1)\times (n-1)$ matrix $B(x)=[b_{ij}(x)]$ such that :
$$\left( \begin{array}{c}
V(x) \cr
V_2(x) \cr
\vdots \cr
V_n(x) \cr
\end{array} \right)= \left( \begin{array}{cc}
1 & 0 \cr 0 & B(x) \cr
\end{array} \right)\left( \begin{array}{c}
V(x) \cr
V_2'(x) \cr
\vdots \cr
V_n'(x) \cr
\end{array} \right).$$
This gives that :
$$\left( \begin{array}{c}
\langle V_2(x), e_1 \rangle \cr
\vdots \cr
\langle V_n(x),e_1 \rangle \cr
\end{array} \right)= B(x) \left( \begin{array}{c}
\langle V_2'(x),e_1 \rangle \cr
\vdots \cr
\langle V_n'(x),e_1 \rangle \cr
\end{array} \right)=B(x) \left( \begin{array}{c}
 -a_2(x) \cr
\vdots \cr
 -a_n(x) \cr
\end{array} \right),
$$
and that $\frac{W_{e_1}}{\vert W_{e_1} \vert}(x)=e_1$. 

If $\frac{W_{e_1}}{\vert W_{e_1} \vert}(x)=e_1$ then $a_1(x)>0$ and $\langle V_i(x),X \rangle =0$ for $i \in \{2,\ldots,n\}$.
This implies that $a_i(x)=0$ for $i \in \{2,\ldots,n\}$ because $B(x)$ is invertible. $\hfill \Box$

Before going further on, we need to carry out some technical computations. Assume that $H=(H_1,\ldots,H_n) : \mathbb{R}^n
\rightarrow \mathbb{R}^n$ is a smooth mapping which does not vanish on a sphere $S_r^{n-1}$. Then we can consider the mapping 
$\frac{H}{\vert H \vert} : S_r^{n-1} \rightarrow S^{n-1}$. Let $P$ be a point in $S_r^{n-1}$ such that $\frac{H}{\vert H \vert}(P)=e_1$. We
can assume that $x_1(P) \not= 0$. If we set $x=(x_1,x')$ where $x'$ belongs to $\mathbb{R}^{n-1}$ then, by the implicit
function theorem, there exists a smooth function $\varphi : \mathbb{R}^{n-1} \rightarrow \mathbb{R}$ such that in the
neighborhood of $P$, $S_r^{n-1}$ is the set of points $(\varphi(x'),x')$. Let us write $\theta(x')=(\varphi(x'),x')$. Let
deg$(\theta,P')$ be the degree of $\theta$ at $P'$ where we write $P=(x_1(P),P')$ ; it is $+1$ if $\theta$ preserves the
orientation and $-1$ otherwise. As explained in [Du1], Lemma 2.2, we have deg$(\theta,P)=\hbox{sign }x_1(P)$. Let $\tilde{H}$
be the mapping defined in the neighborhood of $P'$ by :
$$\tilde{H}(x')=\left( \frac{H_2(\theta(x'))}{\vert H(\theta(x')) \vert},\ldots, 
\frac{H_n(\theta(x'))}{\vert H(\theta(x')) \vert} \right).$$
Since $H_1(P)>0$, we have :
$$\hbox{deg}(\frac{H}{\vert H\vert},P)=\hbox{sign } x_1(P) \ \hbox{deg}(\tilde{H},P').$$
Differentiating the equality :
$$\tilde{H}_i(x')=\frac{H_i(\theta(x'))}{\vert H(\theta(x')) \vert},$$
and using the fact that $H_i(P)=0$, we find that for $(i,j) \in \{2,\ldots,n\}^2$ :
$$\frac{\partial \tilde{H}_i}{\partial x_j}(P')=\frac{1}{\vert H(P) \vert}\left(
\frac{\partial H_i}{\partial x_1}(P)\frac{\partial \varphi}{\partial x_j}(P')+\frac{\partial H_i}{\partial x_j}(P) \right).$$
Finally $P$ is a regular point of $\frac{H}{\vert H \vert} :S_r^{n-1} \rightarrow S^{n-1}$ if and only if :
$$\hbox{det} \left[ \frac{\partial H_i}{\partial x_1}(P)\frac{\partial \varphi}{\partial x_j}(P')+\frac{\partial H_i}{\partial
x_j}(P) \right]_{(i,j) \in \{2,\ldots,n\}^2} \not= 0.$$
In this situation, we have :
$$\hbox{deg}(\frac{H}{\vert H \vert},P)=\hbox{sign } x_1(P) \  
\hbox{det} \left[ \frac{\partial H_i}{\partial x_1}(P)\frac{\partial \varphi}{\partial x_j}(P')+\frac{\partial H_i}{\partial
x_j}(P) \right]_{(i,j) \in \{2,\ldots,n\}^2}.$$

Let us choose $r>0$ small such that $V^{-1}(0) \cap B_r^n =W_{e_1}^{-1}(0) \cap B_r^n =\{0\}$. We know that
$\hbox{Ind}_{PH}(V,0,\mathbb{R}^n)$ is the
topological degree of $\frac{V}{\vert V \vert }: S_r^{n-1} \rightarrow S^{n-1}$ and that $\hbox{Ind}_{PH}(W_{e_1},0,\mathbb{R}^n)$ is the
topological degree of $\frac{W_{e_1}}{\vert W_{e_1} \vert }: S_r^{n-1} \rightarrow S^{n-1}$.

\begin{lemma}
The vector $e_1$ is a regular value of $\frac{V}{\vert V \vert} : S_r^{n-1} \rightarrow S^{n-1}$ if and only if it is a regular value
of $\frac{W_{e_1}}{\vert W_{e_1} \vert} : S_r^{n-1} \rightarrow S^{n-1}$. In this situation, we have :
$$\hbox{\em deg}(\frac{V}{\vert V \vert},P)=(-1)^{n-1} \hbox{\em deg}(\frac{W_{e_1}}{\vert W_{e_1} \vert},P),$$
for all $P$ in $S_r^{n-1}$ such that $\frac{V}{\vert V \vert}(P)=e_1$.
\end{lemma}
{\it Proof.} Let $P$ be a point such that $\frac{V}{\vert V \vert}(P)=\frac{W_{e_1}}{\vert W_{e_1} \vert}(P)=e_1$. With the
notations of the previous lemma, we have for $x$ close to $P$ and  for $i \in \{2,\ldots,n\}$ :
$$\langle V_i(x),e_1 \rangle =-\sum_{k=2}^nb_{ik}(x)a_k(x),$$
hence for $j \in \{1,\ldots,n\}$,
$$\frac{\partial \langle V_i(P),e_1 \rangle}{\partial x_j}=-\sum_{k=2}^n b_{ik}(P) \frac{\partial a_k}{\partial x_j}(P).$$
Applying the above computations to $\frac{V}{\vert V \vert}$ and $\frac{W_{e_1}}{\vert W_{e_1} \vert}$, it is easy to conclude.
$\hfill \Box$

Now we can state the relation between the two indices.
\begin{proposition}
We have :
$$\hbox{\em Ind}_{PH}(V,0,\mathbb{R}^n) =(-1)^{n-1}\hbox{\em Ind}_{PH}(W_{e_1},0,\mathbb{R}^n).$$
\end{proposition}
{\it Proof.} Let us fix $r>0$ such that $V^{-1}(0) \cap B_r^n =W_{e_1}^{-1}(0) \cap B_r^n=\{0\}$. If $e_1$ is a regular value of
$\frac{V}{\vert V \vert}:S_r^{n-1} \rightarrow S^{n-1} $, we combine the two previous lemmas to get the result.

If $e_1$ is not a regular value of $\frac{V}{\vert V \vert}$, we choose a regular value $w$ of $\frac{V}{\vert V \vert}:S_r
\rightarrow S^{n-1}$ very close to $e_1$. There exists a direct orthogonal matrix $A$, close to $I_n$, such that $Aw=e_1$. 
Let $\bar{V}$ be the vector field defined by $\bar{V}=AV$ and, for $i \in
\{2,\ldots,n\}$, let $\bar{V}_i$ be defined by $\bar{V}_i=AV_i$. The vector field $\bar{V}$ satisfies Condition $(P')$ for
$A$ is direct orthogonal and we have :
$$\hbox{Ind}_{PH}(\bar{V},0,\mathbb{R}^n)=\hbox{Ind}_{PH}(V,0,\mathbb{R}^n).$$ Moreover, since $\frac{\bar{V}}{\vert \bar{V} \vert}(x)=e_1$ if and only
if $\frac{V}{\vert V \vert}(x)=w$, $e_1$ is a regular value of $\frac{\bar{V}}{\vert \bar{V} \vert} : S_r \rightarrow S^{n-1}$
and, by the previous case, $\hbox{Ind}_{PH}(\bar{V},0,\mathbb{R}^n)=(-1)^{n-1}\hbox{Ind}_{PH}(\bar{W}_{e_1},0,\mathbb{R}^n)$
where :
$$\bar{W}_{e_1}=\langle \bar{V},e_1 \rangle \frac{\partial}{\partial x_1}+
\langle \bar{V}_2,e_1 \rangle \frac{\partial}{\partial x_2}+\cdots+
\langle \bar{V}_n,e_1 \rangle \frac{\partial}{\partial x_n}.$$
But $\bar{W}_{e_1}$ is equal to the vector field :
$$\langle V,A^te_1 \rangle \frac{\partial}{\partial x_1}+
\langle V_2,A^te_1 \rangle \frac{\partial}{\partial x_2}+\cdots+
\langle V_n,A^te_1 \rangle \frac{\partial}{\partial x_n},$$
whose index at the origin is $\hbox{Ind}_{PH}(W_{e_1},0,\mathbb{R}^n)$ (here $A^t$ is the transpose matrix of $A$). $\hfill \Box$

Let $Z=b_1 \frac{\partial}{\partial x_1}+\cdots+b_n \frac{\partial}{\partial x_n}$ be another smooth vector field defined near the
origin and let $\Gamma$ be the following vector field :
$$\Gamma=\langle V,Z \rangle \frac{\partial}{\partial x_1}+\langle V_2,Z \rangle \frac{\partial}{\partial x_2}+\cdots+\langle V_n,Z
\rangle \frac{\partial}{\partial x_n}.$$
The next theorem relates the indices of $V$, $Z$ and $\Gamma$.
\begin{theorem}
The vector field $\Gamma$ has an isolated zero at the origin if and only if $Z$ has an isolated zero at the origin. In this case, we
have :
$$\hbox{\em Ind}_{PH}(\Gamma,0,\mathbb{R}^n)= \hbox{\em Ind}_{PH}(Z,0,\mathbb{R}^n) +(-1)^{n-1}\hbox{\em Ind}_{PH}(V,0,\mathbb{R}^n).$$
\end{theorem}
{\it Proof.} The equivalence is clear because of Condition $(P')$ and the fact that $V$ has an isolated zero at $0$. To prove the equality, we distinguish two cases. The
first case is when there exists $j \in \{2,\ldots,n\}$ such that $V_j(0) \not= 0$. Let
$\tilde{Z}=\tilde{b}_1\frac{\partial}{\partial x_1}+\cdots+\tilde{b}_n\frac{\partial}{\partial x_n}$ be a small perturbation
of $Z$ such that $\tilde{Z}(0) \notin V_j(0)^\perp$ and the zeros of $\tilde{Z}$ lying close to the origin are non-degenerate.
Let $Q_1,\ldots,Q_s$ be these zeros. Let $\tilde{\Gamma}$ be the vector field defined by :
$$\tilde{\Gamma}=\langle V,\tilde{Z} \rangle\frac{\partial}{\partial x_1}+
\langle V_2,\tilde{Z} \rangle\frac{\partial}{\partial x_2}+\cdots+
\langle V_n,\tilde{Z} \rangle\frac{\partial}{\partial x_n}.$$
The points $Q_1,\ldots,Q_s$ are exactly the zeros of $\tilde{\Gamma}$ near the origin. Let us compare the signs of :
$$\frac{\partial(\langle V,\tilde{Z} \rangle, \langle V_2,\tilde{Z} \rangle, \ldots,\langle V_n,\tilde{Z} \rangle)}
{\partial(x_1,\ldots,x_n)}(Q_j),$$ and :
$$\frac{\partial(\tilde{b}_1,\ldots,\tilde{b}_n)}
{\partial(x_1,\ldots,x_n)}(Q_j),$$
for $j \in \{1,\ldots,s\}$. Since $(V(Q_j),V_2(Q_j),\ldots,V_n(Q_j))$ is a direct basis, the matrix $B(Q_j)$
given by :
$$B(Q_j)=\left( \begin{array}{c}
V(Q_j) \cr V_2(Q_j) \cr \vdots \cr V_n(Q_j) \cr \end{array} \right),$$
is a direct matrix.  A straightforward computation gives that :
$$\displaylines{
\qquad \frac{\partial(\langle V,\tilde{Z} \rangle, \langle V_2,\tilde{Z} \rangle, \ldots,\langle V_n,\tilde{Z} \rangle)}
{\partial(x_1,\ldots,x_n)}(Q_j)= \hfill \cr
 \hbox{det }B(Q_j) \frac{\partial(\langle e_1,\tilde{Z} \rangle, \langle e_2,\tilde{Z} \rangle, 
\ldots,\langle e_n,\tilde{Z} \rangle)}{\partial(x_1,\ldots,x_n)}(Q_j)= \cr
\hfill \hbox{det }B(Q_j)\frac{\partial(\tilde{b}_1,\ldots,\tilde{b}_n)}
{\partial(x_1,\ldots,x_n)}(Q_j). \qquad \cr
}$$
Now $\hbox{Ind}_{PH}(\Gamma,0,\mathbb{R}^n)$ (resp. $\hbox{Ind}_{PH}(Z,0,\mathbb{R}^n)$) is the degree around a small sphere of $\frac{\tilde{\Gamma}}{\vert \tilde{\Gamma}
\vert}$ (resp. $\frac{\tilde{Z}}{\vert \tilde{Z} \vert}$), and the above equality shows that $$\hbox{Ind}_{PH}(\Gamma,0,\mathbb{R}^n)
=\hbox{Ind}_{PH}(Z,0,\mathbb{R}^n) .$$ 
Since $V_j(0) \not= 0$, $\hbox{Ind}_{PH}(V_j,0,\mathbb{R}^n)$ is zero. This index is also the topological degree around a
small sphere $S_r^{n-1}$ of
$\frac{V_j}{\vert V_j \vert}$. But for each point $x$ in $S_r^{n-1}$, $(V(x),V_2(x),\ldots,V_n(x))$ is a direct basis. Hence the
vectors $\frac{V_j}{\vert V_j \vert}(x)$ and $\frac{V}{\vert V \vert}(x)$ are not opposite vectors and the mappings 
$\frac{V_j}{\vert V_j \vert} : S_r^{n-1} \rightarrow S^{n-1} $ and $\frac{V}{\vert V \vert} : S_r^{n-1} \rightarrow S^{n-1}$are homotopic.
Finally $\hbox{Ind}_{PH}(V,0,\mathbb{R}^n)=\hbox{Ind}_{PH}(V_j,0,\mathbb{R}^n)=0$.

Now assume that for all $j \in \{2,\ldots,n\}$, $V_j(0)=0$. Let
$\tilde{Z}=\tilde{b}_1\frac{\partial}{\partial x_1}+\cdots+\tilde{b}_n\frac{\partial}{\partial x_n}$ be a small perturbation
of $Z$ such that $\tilde{Z}(0) \not= 0$ and the zeros of $\tilde{Z}$ lying close to the origin are non-degenerate.
Let $Q_1,\ldots,Q_s$ be these zeros. Let $\tilde{\Gamma}$ be the vector field defined by :
$$\tilde{\Gamma}=\langle V,\tilde{Z} \rangle\frac{\partial}{\partial x_1}+
\langle V_2,\tilde{Z} \rangle\frac{\partial}{\partial x_2}+\cdots+
\langle V_n,\tilde{Z} \rangle\frac{\partial}{\partial x_n}.$$
The zeros of $\tilde{\Gamma}$ are $Q_1,\ldots,Q_s$ and the origin. Furthermore, we have :
$$\hbox{Ind}_{PH}(\Gamma,0,\mathbb{R}^n)=\sum_{j=1}^s \hbox{Ind}_{PH}(\tilde{\Gamma},Q_j,\mathbb{R}^n)
+\hbox{Ind}_{PH}(\tilde{\Gamma},0,\mathbb{R}^n).$$
For the same reasons as in the first case, we have : $$\sum_{j=1}^s \hbox{Ind}_{PH}(\tilde{\Gamma},Q_j,\mathbb{R}^n)
=\hbox{Ind}_{PH}(Z,0,\mathbb{R}^n).$$ Since
$\tilde{Z}(0)\not= 0$, $\hbox{Ind}_{PH}(\tilde{\Gamma},0,\mathbb{R}^n)$ is equal to the index at the origin of the vector
field : 
$$\langle V,X \rangle\frac{\partial}{\partial x_1}+
\langle V_2,X \rangle\frac{\partial}{\partial x_2}+\cdots+
\langle V_n,X \rangle\frac{\partial}{\partial x_n},$$
where $X=\frac{\tilde{Z}(0)}{\vert \tilde{Z}(0) \vert}$. This index is equal to
$(-1)^{n-1}\hbox{Ind}_{PH}(V,0,\mathbb{R}^n)$, by Lemma 6.3 and
Proposition 6.6. $\hfill \Box$

\section{1-forms and hypersurfaces with isolated singularities}
Let $f :(\mathbb{R}^n,0) \rightarrow (\mathbb{R},0)$ be an analytic function defined in the neighborhood of $0$ with an
isolated critical point at the origin. Let $\Omega =a_1 dx_1+\cdots+a_n dx_n$ be a smooth 1-form. In this section,
under some assumptions on $f$ or on $\Omega$, we relate Ind$_{Rad}(\Omega,0,f^{-1}(0))$, Ind$_{Rad}(\Omega,0,\{f\ge 0\})$ 
and Ind$_{Rad}(\Omega,0,\{f\le 0\})$ to usual Poincar\'e-Hopf indices of vector fields.

Let us recall first the following formula due to Khimshiashvili [Kh] and that we will use in our proofs. 
If $\delta$ is a regular value of $f$ such that $0< \vert \delta \vert
\ll r \ll 1$ then we have :
$$\chi(f^{-1}(\delta) \cap B_r^n)=1-\hbox{sign}(-\delta)^n  \hbox{Ind}_{PH}(\nabla f,0,\mathbb{R}^n).$$
Moreover, we also have (see [Du1], Theorem 3.2) :
$$\chi(\{f \ge \delta\} \cap B_r^n)
-\chi(\{f \le \delta\} \cap B_r^n)=\hbox{sign}(-\delta)^{n-1} \hbox{Ind}_{PH}(\nabla f,0,\mathbb{R}^n).$$

As usual, we will work with the coordinate system $(x_1,\ldots,x_n)$. First we assume that the vector field $\nabla f$
satisfies Condition $(P')$ of Section 6 : there exist smooth vector fields $V_2,\ldots,V_n$ such that
$V_2(x),\ldots,V_n(x)$ span $(\nabla f(x))^\perp$, whenever $\nabla f(x)\not= 0$, and such that the orientation of $(\nabla
f(x),V_2(x),\ldots,$ $V_n(x))$ agrees with the orientation of $\mathbb{R}^n$. 

Let $V(\Omega)$ and $W(f,\Omega)$ be the following vector fields :
$$V(\Omega)=a_1 \frac{\partial}{\partial x_1}+\cdots+a_n\frac{\partial}{\partial x_n},$$
$$W(f,\Omega)=f \frac{\partial}{\partial x_1}+\langle V(\Omega),V_2 \rangle \frac{\partial}{\partial x_2}+\cdots+
\langle V(\Omega),V_n \rangle \frac{\partial}{\partial x_n}.$$

\begin{lemma}
The vector field $W(f,\Omega)$ has an isolated zero at the origin if and only if $\Omega$ has an isolated
zero at $0$ on $f^{-1}(0)$.
\end{lemma}
{\it Proof.} The form $\Omega$ has a zero at a point $x$ on $ f^{-1}(0)$ different from the origin
if and only if $f(x)=0$ and $\Omega(x)$ is proportional to $df(x)$. 
This last condition is equivalent to the fact that $\langle
V(x),V_i(x) \rangle$ vanishes for $i\in \{2,\ldots,n\}$. $\hfill \Box$

\begin{theorem}
Assume that $W(f,\Omega)$ has an isolated zero at the origin. Then we have :
$$\hbox{\em Ind}_{Rad}(\Omega,0,f^{-1}(0))= \hbox{\em Ind}_{PH}(\nabla f,0,\mathbb{R}^n) + 
\hbox{\em Ind}_{PH} (W(f,\Omega),0,\mathbb{R}^n) . $$
\end{theorem}
{\it Proof.} Let us fix $r>0$ sufficiently small so that $S_{r'}^{n-1}$ intersects $f^{-1}(0)$ transversally for $0<r'\le r$
and $\Omega$ has no zero on $f^{-1}(0) \setminus\{0\}$ inside $B_r^n$. Let
$\tilde{\Omega}=\tilde{a}_1dx_1+\cdots+\tilde{a}_n dx_n$ be a small perturbation of $\Omega$ such that 
$\tilde{\Omega}$ is a correct and non-degenerate form on $f^{-1}(0) \cap \{r' \le \vert x \vert \le r\}$, for some $r'<r$. Let $\{P_i\}$ be the set of
inward zeros of $\tilde{\Omega}$ on $f^{-1}(0) \cap \{r' \le \vert x \vert \le r\}$ lying in $S_r^{n-1}$. 
Using the same method as in Lemma 5.2, we
can prove that :
$$\hbox{Ind}_{Rad}(\Omega,0,f^{-1}(0))=1-\sum_i 
\hbox{Ind}_{PH}(\tilde{\Omega},P_i,S_r^{n-1}\cap f^{-1}(0)) .$$
We can also assume that if $\delta \not= 0$ is small enough then $\tilde{\Omega}$ is
correct and non-degenerate on $f^{-1}(\delta) \cap B_r^n$. Let us denote by $Q_1,\ldots,Q_s$ its singular points not lying in
 $f^{-1}(\delta) \cap S_r^{n-1}$. By the
Poincar\'e-Hopf theorem, we have :
$$\chi(f^{-1}(\delta) \cap B_r^n)=\sum_{i=1}^s \hbox{Ind}_{PH}(\tilde{\Omega},Q_i,f^{-1}(\delta)) + 
1-\hbox{Ind}_{Rad}(\Omega,0,f^{-1}(0)).$$
So we have to relate the sum of indices in the right-hand side of this equality to the index of $W(f,\Omega)$. Let us fix $i$
in $\{1,\ldots,s\}$ and let us set $Q=Q_i$ for convenience. Since $\delta$ is a regular value of $f$, there exists $j$ such
that $\frac{\partial f}{\partial x_j}(Q) \not= 0$. Assume that $j=1$. By Lemma 2.1, we have :
$$\hbox{Ind}_{PH}(\tilde{\Omega},Q,f^{-1}(\delta))=\hbox{sign} \left( (-1)^{n-1} \frac{\partial f}{\partial x_1}(Q)^n 
\frac{\partial(f,\tilde{m}_1,\ldots,\tilde{m}_n)}{\partial(x_1,\ldots,x_n)}(Q) \right),$$
where for $j \ge 2$,
$$\tilde{m}_j =\left\vert \begin{array}{cc}
\tilde{a}_1 & \tilde{a}_j \cr
\frac{\partial f}{\partial x_1} & \frac{\partial f}{\partial x_j} \cr
\end{array} \right\vert.$$
A computation, similar to the one done in [Du3,Lemma 2.5] in the case of the differential of a function, gives :
$$\displaylines{
\qquad \hbox{sign} \left(  
\frac{\partial(f-\delta,\tilde{m}_1,\ldots,\tilde{m}_n)}{\partial(x_1,\ldots,x_n)}(Q) \right)= \hfill \cr
\hfill \hbox{sign} \left( (-1)^{n-1} \frac{\partial f}{\partial x_1}(Q)^n 
\frac{\partial(f-\delta,\langle \tilde{V},V_2 \rangle,\ldots,\langle \tilde{V},V_n \rangle)}{\partial(x_1,\ldots,x_n)}(Q) \right)
,\qquad
}$$
where $\tilde{V}=\tilde{a}_1 \frac{\partial}{\partial x_1}+\cdots+\tilde{a}_n \frac{\partial}{\partial x_n}$. This proves
that :
$$\hbox{Ind}_{PH}(\tilde{\Omega},Q,f^{-1}(\delta))=\hbox{sign} \left(
\frac{\partial(f,\langle \tilde{V},V_2 \rangle,\ldots,\langle \tilde{V},V_n \rangle)}{\partial(x_1,\ldots,x_n)}(Q)
\right).$$
Summing over all the points $Q_i$, we find that $\sum_{i=1}^s \hbox{Ind}_{PH}(\tilde{\Omega},Q_i,f^{-1}(\delta))$ is equal to the
degree of the mapping $\frac{\tilde{W}}{\vert \tilde{W} \vert} : S_r^{n-1} \rightarrow S^{n-1}$, where
$\tilde{W}=f\frac{\partial}{\partial x_1}+\langle \tilde{V},V_2 \rangle \frac{\partial}{\partial x_2}+\cdots+
\langle \tilde{V},V_n \rangle \frac{\partial}{\partial x_n},$ which is equal to Ind$_{PH}(W(f,\Omega),0,\mathbb{R}^n)$. Hence :
$$\hbox{Ind}_{Rad}(\Omega,0,f^{-1}(0))=1-\chi(f^{-1}(\delta) \cap B_r^n)+\hbox{Ind}_{PH}(W(f,\Omega),0,\mathbb{R}^n).$$
To end the proof, we apply Khimshiashvili's formula. If $n$ is even, $\chi(f^{-1}(\delta) \cap B_r^n)=1-\hbox{Ind}_{PH}
(\nabla f,0,\mathbb{R}^n)$.
If $n$ is odd, Ind$_{PH}(\nabla f,0,\mathbb{R}^n)=0$ as recalled in Section 6 and $\chi(f^{-1}(\delta) \cap B_r^n)=1$. $\hfill \Box$

We can apply Theorem 7.2 to the differential of an analytic function and recover the results of Theorem 2.1 in [Du3].
\begin{corollary}
Let $g:(\mathbb{R}^n,0) \rightarrow (\mathbb{R},0)$ be an analytic function defined in the neighborhood of the origin such that
$g(0)=0$. Let us assume that $g_{\vert f^{-1}(0) \setminus \{0\}}$ has no critical point in the neighborhood of the origin.
Then the vector field $W(f,dg)$ has an isolated zero at the origin. If $n$ is even, we have:
$$\displaylines{
\qquad \chi(f^{-1}(0) \cap g^{-1}(\delta) \cap B_r^n)= 1-\hbox{\em Ind}_{PH}(\nabla f,0,\mathbb{R}^n) + \hfill \cr
\hfill \hbox{\em sign}(\delta) \hbox{\em Ind}_{PH}(W(f,dg),0,\mathbb{R}^n). \qquad \cr
}$$
If $n$ is odd, we have :
$$\chi(f^{-1}(0) \cap g^{-1}(\delta) \cap B_r^n)= 1-\hbox{\em Ind}_{PH}(W(f,dg),0,\mathbb{R}^n) .$$
\end{corollary}
{\it Proof.} Combine Theorem 7.2 and Theorem 2 in [EG5]. $\hfill \Box$

Now let us study Ind$_{Rad}(\Omega,0,\{f\ge 0\})$ and Ind$_{Rad}(\Omega,0,\{f\le 0\})$. Let $Y(f,\Omega)$ and 
$\Gamma(f,\Omega)$ be the following vector fields :
$$Y(f,\Omega)=f\langle \nabla f,V(\Omega) \rangle \frac{\partial}{\partial x_1}+ \langle V(\Omega),V_2 \rangle \frac{\partial}{\partial
x_2}+\cdots+\langle V(\Omega),V_n \rangle \frac{\partial}{\partial x_n},$$
$$\Gamma(f,\Omega)=\langle \nabla f,V(\Omega) \rangle \frac{\partial}{\partial x_1}+ \langle V(\Omega),V_2 \rangle \frac{\partial}{\partial
x_2}+\cdots+\langle V(\Omega),V_n \rangle \frac{\partial}{\partial x_n}.$$  
\begin{lemma}
The vector field $Y(f,\Omega)$ has an isolated zero at the origin if and only if the vector fields $V(\Omega)$ and
$W(f,\Omega)$ have an isolated zero at the origin.
\end{lemma}
{\it Proof.} It has an isolated zero at the origin if and only if $W(f,\Omega)$ and $\Gamma(f,\Omega)$ have an isolated zero at
the origin. It is enough to apply the first assertion of Theorem 6.7. $\hfill \Box$

\begin{lemma}
The form $\Omega$ has an isolated zero at the origin on $ \{f \ge 0\}$ and $\{ f \le 0\}$ if and only if
$Y(f,\Omega)$ has an isolated zero at the origin.
\end{lemma}
{\it Proof.} This is easy using the previous lemma and proceeding as in Lemma 7.1. $\hfill \Box$

\begin{theorem}
Assume that $Y(f,\Omega)$ has an isolated zero at the origin. Then we have :
$$\displaylines{
\quad \hbox{\em Ind}_{Rad}(\Omega,0,\{f \ge 0\})=\frac{1}{2} \Big[ \hbox{\em Ind}_{PH}(V(\Omega),0,\mathbb{R}^n)+ 
\hfill \cr
\hfill  \hbox{\em Ind}_{PH}(W(f,\Omega),0,\mathbb{R}^n)+ 
\hbox{\em Ind}_{PH} (\nabla f,0,\mathbb{R}^n)+
\hbox{\em Ind}_{PH} (Y(f,\Omega),0,\mathbb{R}^n) \Big], \quad \cr
}$$
$$\displaylines{
\quad \hbox{\em Ind}_{Rad}(\Omega,0,\{f \le 0\})=\frac{1}{2} \Big[ \hbox{\em Ind}_{PH}(V(\Omega),0,\mathbb{R}^n)+ 
\hfill \cr
\hfill  \hbox{\em Ind}_{PH}(W(f,\Omega),0,\mathbb{R}^n)+
\hbox{\em Ind}_{PH} (\nabla f,0,\mathbb{R}^n)-
\hbox{\em Ind}_{PH} (Y(f,\Omega),0,\mathbb{R}^n) \Big]. \quad \cr
}$$
\end{theorem}
{\it Proof.} Let us fix $r>0$ sufficiently small so that $S_{r'}^{n-1}$ intersects $f^{-1}(0)$ transversally for $0<r'\le r$,
$\Omega_{\vert f^{-1}(0) \setminus \{0\}}$ has no zero inside $B_r^n$ and $\Omega$ has no zero on $B_r^n$ except $0$. 

Let $\tilde{\Omega}=\tilde{a}_1dx_1+\cdots+\tilde{a}_ndx_n$ be a  small perturbation of $\Omega$ such that
$\tilde{\Omega}$ is correct and non-degenerate on $\{ f\ge 0\} \cap \{ r' \le \vert x \vert \le r\}$ and 
on $\{ f\le 0\} \cap \{ r' \le \vert x \vert \le r\}$. As above, we denote by $\tilde{V}$ the vector field dual to
$\tilde{\Omega}$. 
Let $\{R_k\}$ (resp. $\{S_l\}$) be the set of inward zeros of $\tilde{\Omega}$ on $\{ f\ge 0\} \cap \{ r' \le \vert x
\vert \le r\}$ (resp. $\{ f\le 0\} \cap \{ r' \le \vert x \vert \le r\}$) lying on $S_r^{n-1}$. 
Using the same method as in
Lemma 5.2, we can prove that :
$$\hbox{Ind}_{Rad}(\Omega,0,\{f \ge 0\})=1-\sum_k 
\hbox{Ind}_{PH} (\tilde{\Omega},R_k,\{f \ge 0 \} \cap S_r^{n-1}),$$
$$\hbox{Ind}_{Rad} (\Omega,0,\{f \le 0\})=1-\sum_l 
\hbox{Ind}_{PH} (\tilde{\Omega},S_l,\{f \le 0 \} \cap S_r^{n-1}).$$
We can also assume that if $\delta \not=0$ is small enough then $\tilde{\Omega}_{\vert \{f \ge \delta\} \cap B_r^n}$ and 
$\tilde{\Omega}_{\vert \{f \le \delta\} \cap B_r^n}$ are correct and non-degenerate and that the zeros of
$\tilde{\Omega}$ lie in $\{ \vert f \vert < \delta \} \cap \mathring{B_r^n}$, where $\mathring{B_r^n}$ is the interior of 
$B_r^n$. Let us denote by $P_1,\ldots,P_s$ the singular
points of $\tilde{\Omega}$ lying in $\mathring{B_r^n}$ and by $Q_1,\ldots,Q_t$ the singular points of
$\tilde{\Omega}_{\vert f^{-1}(\delta) \cap \mathring{B_r^n}}$. By the Poincar\'e-Hopf theorem, we have :
$$\displaylines{
\qquad \chi(\{ f \ge \delta \} \cap B_r^n)= \sum_{j \vert \langle \nabla f(Q_j),\tilde{V}(Q_j) \rangle >0} 
\hbox{Ind}_{PH} (\tilde{\Omega},Q_j,f^{-1}(\delta)) +\hfill \cr
\sum_{i \vert f(P_i) > \delta}
\hbox{Ind}_{PH} (\tilde{\Omega},P_i,\mathbb{R}^n) +1 -\hbox{Ind}_{Rad}(\Omega,0,\{f \ge 0 \}),
}$$
$$\displaylines{
\qquad \chi(\{ f \le \delta \} \cap B_r^n)= \sum_{j \vert \langle \nabla f(Q_j),\tilde{V}(Q_j) \rangle <0} 
\hbox{Ind}_{PH}(\tilde{\Omega},Q_j,f^{-1}(\delta)) + \hfill \cr
\sum_{i \vert f(P_i) < \delta}
\hbox{Ind}_{PH} (\tilde{\Omega},P_i,\mathbb{R}^n) +1 -\hbox{Ind}_{Rad}(\Omega,0,\{f \le 0 \}).
}$$
Summing these two equalities and using the Mayer-Vietoris sequence, we obtain :
$$\displaylines{
\qquad \hbox{Ind}_{Rad}(\Omega,0,\{f \ge 0 \})+\hbox{Ind}_{Rad}(\Omega,0,\{f \le 0 \})= \hfill \cr
\hfill \sum_{j } \hbox{Ind}_{PH}(\tilde{\Omega},Q_j,f^{-1}(\delta)) + \sum_{i }
\hbox{Ind}_{PH}(\tilde{\Omega},P_i,\mathbb{R}^n)+1-\chi(f^{-1}(\delta) \cap B_r^n). \qquad
\cr}$$
As explained in Theorem 7.2 : 
$$\sum_{j } \hbox{Ind}_{PH}(\tilde{\Omega},Q_j,f^{-1}(\delta))=\hbox{Ind}_{PH}(W(f,\Omega),0,\mathbb{R}^n),$$ and 
$\sum_{i }\hbox{Ind}_{PH}(\tilde{\Omega},P_i,\mathbb{R}^n)$ is clearly equal to Ind$_{PH}(V(\Omega),0,\mathbb{R}^n)$. Finally, we have :
$$\displaylines{
\qquad \hbox{Ind}_{Rad}(\Omega,0,\{f \ge 0 \})+\hbox{Ind}_{Rad}(\Omega,0,\{f \le 0 \})=\ \hfill \cr
\hfill \hbox{Ind}_{PH}(V(\Omega),0,\mathbb{R}^n)+\hbox{Ind}_{PH}(W(f,\Omega),0,\mathbb{R}^n)
+\hbox{Ind}_{PH}(\nabla f,0,\mathbb{R}^n).
}$$
Making the difference of the two above equalities leads to :
$$\displaylines{
\qquad \hbox{Ind}_{Rad}(\Omega,0,\{f \ge 0 \})-\hbox{Ind}_{Rad}(\Omega,0,\{f \le 0 \})=\hfill \cr
\hfill \sum_{j } \hbox{sign} \langle \nabla f(Q_j),\tilde{V}(Q_j) \rangle 
\hbox{Ind}_{PH} (\tilde{\Omega},Q_j,f^{-1}(\delta)) + \qquad \qquad \qquad \qquad \cr
\hfill \sum_{i } \hbox{sign}(f(P_i)-\delta)
\hbox{Ind}_{PH} (\tilde{\Omega},P_i,\mathbb{R}^n)- \qquad \qquad \qquad \cr
\hfill \left[ \chi(\{f \ge \delta\} \cap B_r^n)
-\chi(\{f \le \delta\} \cap B_r^n) \right]. \qquad \cr
}$$
Since sign$(f(P_i)-\delta)=\hbox{sign}(-\delta)$ for all $i\in \{1,\ldots,s\}$ and :
$$\chi(\{f \ge \delta \} \cap B_r^n)
-\chi(\{f \le \delta \}\cap B_r^n)=\hbox{sign}(-\delta)^{n-1} \hbox{Ind}_{PH}(\nabla f,0,\mathbb{R}^n),$$
we have :
$$\displaylines{
\qquad \hbox{Ind}_{Rad}(\Omega,0,\{f \ge 0 \})-\hbox{Ind}_{Rad}(\Omega,0,\{f \le 0 \})=\hfill \cr
\hfill \sum_{j } \hbox{sign} \langle \nabla f(Q_j),\tilde{V}(Q_j) \rangle 
\hbox{Ind}_{PH} (\tilde{\Omega},Q_j,f^{-1}(\delta)) + \qquad \qquad \qquad \cr
\hfill \hbox{sign}(-\delta) \hbox{Ind}_{PH}(V(\Omega),0,\mathbb{R}) 
-\hbox{sign}(-\delta)^{n-1} \hbox{Ind}_{PH}(\nabla f,0,\mathbb{R}^n). \qquad \cr
}$$
Let $\tilde{Y}$ and $\tilde{\Gamma}$ be the following vector fields :
$$\tilde{Y}=(f-\delta)\langle \nabla f,\tilde{V} \rangle \frac{\partial}{\partial x_1}+ \langle \tilde{V},V_2 \rangle \frac{\partial}{\partial
x_2}+\cdots+\langle \tilde{V},V_n \rangle \frac{\partial}{\partial x_n},$$
$$\tilde{\Gamma}=\langle \nabla f,\tilde{V} \rangle \frac{\partial}{\partial x_1}+ \langle \tilde{V},V_2 \rangle \frac{\partial}{\partial
x_2}+\cdots+\langle \tilde{V},V_n \rangle \frac{\partial}{\partial x_n}.$$  
The zeros of $\tilde{Y}$ are the points $Q_j$'s, $P_i$'s and possibly the origin (see Theorem 6.7). It is easy to see that the
$Q_j$'s are non-degenerate and that :
$$\hbox{Ind}_{PH} (\tilde{Y},Q_j,\mathbb{R}^n)=\hbox{sign} \langle \nabla f (Q_j),\tilde{V}(Q_j) \rangle
\hbox{Ind}_{PH}(\tilde{\Omega},Q_j,f^{-1}(\delta)).$$
By the position of the points $P_i$, we have :
$$\displaylines{
 \hbox{Ind}_{PH}(Y(f,\Omega),0,\mathbb{R}^n)=
\sum_j \hbox{sign}(\langle \nabla f(Q_j),\tilde{V}(Q_j) \rangle \hbox{Ind}_{PH}(\tilde{\Omega},Q_j,f^{-1}(\delta))+ \hfill
\cr
\hfill \hbox{sign}(-\delta) \left[ \sum_i \hbox{Ind}_{PH}(\tilde{\Gamma},P_i,\mathbb{R}^n)
+\hbox{Ind}_{PH}(\tilde{\Gamma},0,\mathbb{R}^n) \right]= \qquad \qquad \cr
\qquad  \sum_j \hbox{sign}(\langle \nabla f(Q_j),\tilde{V}(Q_j) \rangle
\hbox{Ind}_{PH}(\tilde{\Omega},Q_j,f^{-1}(\delta))+ \hfill \cr
\hfill \hbox{sign}(-\delta)\hbox{Ind}_{PH}(\Gamma(f,\Omega),0,\mathbb{R}^n) = \qquad \cr
\qquad \sum_j \hbox{sign}(\langle \nabla f(Q_j),\tilde{V}(Q_j) \rangle) \hbox{Ind}_{PH}(\tilde{\Omega},Q_j,f^{-1}(\delta))+
\hfill \cr
\hfill \hbox{sign}(-\delta) \hbox{Ind}_{PH}(V(\Omega),0,\mathbb{R}^n) +
(-1)^{n-1}\hbox{sign}(-\delta) \hbox{Ind}_{PH}(\nabla f,0,\mathbb{R}^n). \qquad \cr
}$$
Combining all these equalities and using the fact that Ind$_{PH}(\nabla f,0,\mathbb{R}^n) =0$ if $n$ is odd, we find that :
$$\displaylines{
\qquad \hbox{Ind}_{Rad}(\Omega,0,\{f \ge 0 \})-\hbox{Ind}_{Rad}(\Omega,0,\{f \le 0 \})= \hfill \cr
\hfill \hbox{Ind}_{PH}(Y(f,\Omega),0,\mathbb{R}^n) -\hbox{sign}(-\delta)^{n-1}
\hbox{Ind}_{PH}(\nabla f,0,\mathbb{R}^n) -\qquad \qquad \cr
\hfill (-1)^{n-1}\hbox{sign}(-\delta) \hbox{Ind}_{PH}(\nabla f,0,\mathbb{R}^n)= \qquad \cr
\qquad \hbox{Ind}_{PH}(Y(f,\Omega),0,\mathbb{R}^n) -(-1)^{n-1} \hbox{sign}(\delta)^{n-1}
\hbox{Ind}_{PH}(\nabla f,0,\mathbb{R}^n) -\qquad \qquad \cr
\hfill (-1)^{n-1}\hbox{sign}(-\delta) \hbox{Ind}_{PH}(\nabla f,0,\mathbb{R}^n)= \qquad \cr
\qquad \hbox{Ind}_{PH}(Y(f,\Omega),0,\mathbb{R}^n) - \hfill \cr
\qquad \qquad \qquad  (-1)^{n-1} \left[ \hbox{sign}(\delta)^{n-1}
-\hbox{sign}(\delta) \right] \hbox{Ind}_{PH}(\nabla f,0,\mathbb{R}^n)= \hfill \cr
 \hfill \hbox{Ind}_{PH}(Y(f,\Omega),0,\mathbb{R}^n). \qquad \cr
}$$
$\hfill \Box$

\begin{corollary}
Let $g :(\mathbb{R}^n,0) \rightarrow (\mathbb{R},0)$ be an analytic function defined in the neighborhood of the origin such
that $g(0)=0$. Let us assume that $g$ has no critical point on $\{f \ge 0\}$ and $\{f \le 0\}$ in the neighborhood of 
the origin. Then the vector fields $\nabla g$, $W(f,dg)$ and $Y(f,dg)$ have an isolated zero at the origin and if $n$ is even,
we have :
$$\displaylines{
  \chi\left( g^{-1}(\delta) \cap \{ f \ge 0 \} \cap B_r^n \right)= 1-\frac{1}{2}\Big[ 
\hbox{\em Ind}_{PH}(\nabla g,0,\mathbb{R}^n) +\hbox{\em Ind}_{PH}(\nabla f,0,\mathbb{R}^n) +  \hfill \cr 
\hfill  \hbox{\em Ind}_{PH}
(Y(f,dg),0,\mathbb{R}^n) \Big] + \frac{1}{2}\hbox{\em sign}(\delta) \hbox{\em Ind}_{PH}(W(f,dg),0,\mathbb{R}^n) ,
\quad \cr
  \chi\left( g^{-1}(\delta) \cap \{ f \le 0 \} \cap B_r^n \right)= 1-\frac{1}{2}\Big[ 
\hbox{\em Ind}_{PH}(\nabla g,0,\mathbb{R}^n) +\hbox{\em Ind}_{PH}(\nabla f,0,\mathbb{R}^n) -  \hfill \cr
\hfill  \hbox{\em Ind}_{PH}
(Y(f,dg),0,\mathbb{R}^n) \Big] - 
\frac{1}{2}\hbox{\em sign}(\delta) \hbox{\em Ind}_{PH}(W(f,dg),0,\mathbb{R}^n) . \quad \cr
}$$
If $n$ is odd, we have :
$$\displaylines{
\chi\left( g^{-1}(\delta) \cap \{ f \ge 0 \} \cap B_r^n \right)= 1+\frac{1}{2}\hbox{\em sign}(\delta) \Big[ 
\hbox{\em Ind}_{PH}(\nabla g,0,\mathbb{R}^n) +  \hfill \cr
\hfill  \hbox{\em Ind}_{PH}
(Y(f,dg),0,\mathbb{R}^n) \Big] - \frac{1}{2}\hbox{\em Ind}_{PH}(W(f,dg),0,\mathbb{R}^n) , \quad \cr
\chi\left( g^{-1}(\delta) \cap \{ f \le 0 \} \cap B_r^n \right)= 1+\frac{1}{2}\hbox{\em sign}(\delta) \Big[
\hbox{\em Ind}_{PH}(\nabla g,0,\mathbb{R}^n)-  \hfill \cr
\hfill  \hbox{\em Ind}_{PH}
(Y(f,dg),0,\mathbb{R}^n) \Big] - 
 \frac{1}{2}\hbox{\em Ind}_{PH}(W(f,dg),0,\mathbb{R}^n) . \quad \cr
}$$
\end{corollary}
{\it Proof.} Use Theorem 2 in [EG5]. $\hfill \Box$

Now we assume that the vector field $V(\Omega)=a_1 \frac{\partial}{\partial x_1}+\cdots+a_n \frac{\partial}{\partial x_n}$ satisfies
Condition $(P')$ of Section 6 : there exist smooth vector fields $V_2,\ldots,V_n$ in $\mathbb{R}^n$ such that
$V_2(x),\ldots,V_n(x)$ span $[V(\Omega)(x)]^\perp$ whenever $V(\Omega)(x)\not= 0$ and such that $(V(\Omega)(x),V_2(x),\ldots,V_n(x))$ is a direct basis. We
also assume that $\Omega$ (and $V(\Omega)$) has an isolated zero at the origin. Let us consider the following vector fields : 
$$W(f,\Omega)= f\frac{\partial}{\partial x_1}+\langle \nabla f,V_2 \rangle \frac{\partial}{\partial x_2}+\cdots+ 
\langle \nabla f ,V_n \rangle \frac{\partial}{\partial x_n},$$
$$\Gamma(f,\Omega)= \langle \nabla f, V(\Omega) \rangle \frac{\partial}{\partial x_1}+\langle \nabla f,V_2 \rangle \frac{\partial}{\partial x_2}+\cdots+ 
\langle \nabla f ,V_n \rangle \frac{\partial}{\partial x_n},$$
$$Y(f,\Omega)= f\langle \nabla f, V(\Omega) \rangle \frac{\partial}{\partial x_1}+\langle \nabla f,V_2 \rangle \frac{\partial}{\partial x_2}+\cdots+ 
\langle \nabla f ,V_n \rangle \frac{\partial}{\partial x_n}.$$

\begin{lemma}
The vector field $W(f,\Omega)$ has an isolated at $0$ if and only if $\Omega$ has an isolated zero at $0$ on $f^{-1}(0)$.
\end{lemma}
{\it Proof.} See Lemma 7.1. $\hfill \Box$

\begin{lemma}
We can choose $\delta$ small enough and we can perturb $f$ into $\tilde{f}$ in such a way that $\Omega$ has only non-degenerate 
zeros on $f^{-1}(\delta) \cap B_r^n$.
\end{lemma}
{\it Proof.} Let $(x,t)=(x_1,\ldots,x_n,t_1,\ldots,t_n)$ be a coordinate system of $\mathbb{R}^{2n}$ and let :
$$\bar{f}(x,t)=f(x)+\sum_{i=1}^n t_i x_i.$$
For $(i,j) \in \{1,\ldots,n\}^2$, we define $M_{ij}(x,t)$ by :
$$M_{ij}(x,t)= \left\vert \begin{array}{cc}
a_i(x) & a_j(x) \cr
\frac{\partial \bar{f}}{\partial x_i}(x,t) & \frac{\partial \bar{f}}{\partial x_j}(x,t) \cr
\end{array} \right\vert.$$
Notice that :
$$M_{ij}(x,t)= \left\vert \begin{array}{cc}
a_i(x) & a_j(x) \cr
\frac{\partial f}{\partial x_i}(x,t) & \frac{\partial f}{\partial x_j}(x,t) \cr
\end{array} \right\vert +a_it_j-t_ia_j.$$
Let $N$ be defined by :
$$N=\left\{ (x,t) \in \mathbb{R}^{2n} \ \vert \ M_{ij}(x,t)=0 \hbox{ for }(i,j) \in \{1,\ldots,n\}^2 \right\}.$$
At a point $p\not= 0$, $\Omega$ does not vanish, so there exists $i\in \{1,\ldots,n\}$ such that $a_i(p) \not= 0$. This
implies that $N\setminus \{(0,t) \ \vert \ t \in \mathbb{R}^{n} \}$ is a smooth manifold of dimension $n+1$ (or empty).
Actually if $(p,t)$ belongs to $N\setminus \{(0,t) \ \vert \ t \in \mathbb{R}^{n} \}$ then we can assume that $a_1(p)
\not=0$. In this case around $(p,t)$, $N$ is defined by the vanishing of $M_{12},\ldots,M_{1n}$ and the gradient
vectors of these functions are linearly independent. Let $\pi$ be the following mapping :
$$\begin{array}{ccccc}
\pi & : & N\setminus \{(0,t) \ \vert \ t \in \mathbb{R}^{n} \} & \rightarrow & \mathbb{R}^{n+1} \cr
  &   & (x,t) & \mapsto & (\bar{f}(x,t),t) .\cr
\end{array}$$
By the Bertini-Sard theorem, we can choose $(\delta,s)$ close to $0$ in $\mathbb{R}^{n+1}$ such that $\pi$ is regular at
each point in $\pi^{-1}(\delta,s)$. If we denote by $\tilde{f}$ the function defined by $\tilde{f}(x)=f(x,s)$, this means
that $\Omega$ admits on $\tilde{f}^{-1}(\delta)$ only non-degenerate zeros in the neighborhood of the origin. $\hfill
\Box$

\begin{theorem}
Assume that $Y(f,\Omega)$ has an isolated zero at the origin. Then $W(f,\Omega)$ and $\Gamma(f,\Omega)$ also have an
isolated zero at the origin. Furthermore, we have :
$$\displaylines{
\qquad \hbox{if } n \hbox{ is even, } \hbox{\em Ind}_{Rad}(\Omega,0,f^{-1}(0))=\hbox{\em Ind}_{PH}(\nabla f,0,\mathbb{R}^n) 
-\hfill \cr
\hfill \hbox{\em Ind}_{PH}(W(f,\Omega),0,\mathbb{R}^n), \qquad \cr
\qquad \hbox{if } n \hbox{ is odd, } \hbox{\em Ind}_{Rad}(\Omega,0,f^{-1}(0))=\hbox{\em Ind}_{PH}(Y(f,\Omega),0,\mathbb{R}^n). 
\hfill \cr
}$$
\end{theorem}
{\it Proof.} We proceed as in Theorem 7.2. Let us fix $r>0$ sufficiently small so that $S_{r'}^{n-1}$ intersects
$f^{-1}(0)$ transversally for $0<r'\le r$ and $\Omega$ has no zero on $f^{-1}(0) \setminus \{0\}$ inside $B_r^n$. 
By the previous lemma, we can assume that $\Omega$ is correct and non-degenerate on 
$f^{-1}(\delta) \cap B_r^n$.
Morevover, we can assume also that the zeros of $\Omega$ on $f^{-1}(\delta) \cap B_r^n$ lie in $B_{\frac{r}{2}}^n$. Let
us denote them by $Q_1,\ldots,Q_s$. Now we can move $\Omega$ a little in the neighborhood of $f^{-1}(0) \cap S_r^{n-1}$ in
such a way that $\Omega$ is correct on $f^{-1}(0) \cap \{ \frac{3}{4}r \le \vert x \vert \le r \}$ and that no new zeros of 
$\Omega$ are created. As in the proof of Theorem 7.2, we have:
$$\chi(f^{-1}(\delta) \cap B_r^n)= \sum_{i=1}^s \hbox{Ind}_{PH}(\Omega,Q_i,f^{-1}(\delta)) +
1-\hbox{Ind}_{Rad}(\Omega,0,f^{-1}(0)).$$
Let us choose $i \in \{1,\ldots,s\}$ and let us put $Q=Q_i$. Since $\Omega(Q) \not= 0$, there exists $j$ such that
$a_j(Q)\not= 0$. Assume that $j=1$. This implies that $\frac{\partial f}{\partial x_1}(Q) \not= 0$ and by Lemma 2.1, we have
:
$$\hbox{Ind}_{PH}(\Omega,Q,f^{-1}(\delta))= \hbox{sign} \left( (-1)^{n-1} \frac{\partial f}{\partial x_1}(Q)^n
\frac{\partial(f-\delta,m_2,\ldots,m_n)}{\partial(x_1,\ldots,x_n)}(Q) \right),$$
where $m_j=\left\vert \begin{array}{cc}
a_1 & a_j \cr \frac{\partial f}{\partial x_1} & \frac{\partial f}{\partial x_j} \cr
\end{array} \right\vert$. Using the same method as the one used in [Du3], Lemma 2.5 and 2.13 and in Theorem 7.2, we find
that :
$$\displaylines{
\qquad \hbox{sign}\left( \frac{\partial(f-\delta,m_2,\ldots,m_n)}{\partial(x_1,\ldots,x_n)}(Q) \right)= \hfill \cr
\hfill \hbox{sign} \left( a_1(Q)^{n-2} \frac{\partial(f-\delta,\langle \nabla f , V_2 \rangle,\ldots,\langle \nabla f , V_n
\rangle)}{\partial(x_1,\ldots,x_n)}(Q) \right). \qquad \cr
}$$ 
This gives that :
$$\displaylines{
\qquad \hbox{Ind}_{PH}(\Omega,Q,f^{-1}(\delta))= \hfill \cr
\hfill (-1)^{n-1} \hbox{sign} \left( \langle \nabla f(Q),V(Q) \rangle^n
\frac{\partial(f-\delta,\langle \nabla f,V_2 \rangle,\ldots,\langle \nabla f ,V_n \rangle)}{\partial(x_1,\ldots,x_n)}(Q)
\right). \qquad \cr
}$$
When $n$ is even, the proof is the same as in Theorem 7.2. When $n$ is odd, we can relate $\sum_{i=1}^s
\hbox{Ind}_{PH} (\Omega,Q_i,f^{-1}(\delta))$ to Ind$_{PH}(Y(f,\Omega),0,\mathbb{R}^n)$. More precisely, 
as in Theorem 7.6, we have :
$$\displaylines{
\qquad \hbox{Ind}_{PH}(Y(f,\Omega),0,\mathbb{R}^n)= \sum_{i=1}^s \hbox{Ind}_{PH}(\Omega,Q_i,f^{-1}(\delta)) +\hfill \cr
\hfill \hbox{sign}(-\delta)
\hbox{Ind}_{PH}(\Gamma(f,\Omega),0,\mathbb{R}^n), \qquad \cr
}$$
and, by Theorem 6.7 :
$$\displaylines{
\qquad \hbox{Ind}_{PH}(Y(f,\Omega),0,\mathbb{R}^n)= \hfill \cr
\hfill \sum_{i=1}^s \hbox{Ind}_{PH}(\Omega,Q_i,f^{-1}(\delta)) +\hbox{sign}(-\delta) 
\hbox{Ind}_{PH}(\nabla f,0,\mathbb{R}^n). \qquad \cr
}$$
Collecting these informations and using Khimshiashvili's formula, we get :
$$\hbox{Ind}_{Rad}(\Omega,0,f^{-1}(0)) =\hbox{Ind}_{PH}(Y(f,\Omega),0,\mathbb{R}^n).$$
$\hfill \Box$

\begin{corollary}
Let $g :(\mathbb{R}^n,0) \rightarrow (\mathbb{R},0)$ be an analytic function defined in the neighborhood of the origin with
$g(0)=0$. Let us assume that $\nabla g$ satisfies Condition $(P')$ and that $Y(f,dg)$ has an isolated zero at the
origin. Then, if $n$ is even, we have :
$$\displaylines{
\qquad \chi(f^{-1}(0) \cap g^{-1}(\delta) \cap B_r^n)=
1-\hbox{\em Ind}_{PH}(\nabla f,0,\mathbb{R}^n) -\hfill \cr
\hfill \hbox{\em sign}(\delta) \hbox{\em Ind}_{PH}(W(f,dg),0,\mathbb{R}^n). 
\qquad\cr}$$
If $n$ is odd, we have :
$$\displaylines{
\qquad \chi(f^{-1}(0) \cap g^{-1}(\delta) \cap B_r^n)=
1-\hbox{\em Ind}_{PH}(Y(f,dg),0,\mathbb{R}^n). \cr
}$$
$\hfill \Box$
\end{corollary}

Let us study Ind$_{Rad} (\Omega,0,\{f \ge 0\})$ and Ind$_{Rad}(\Omega,0,\{f \le 0 \})$. 

\begin{lemma}
The vector field $Y(f,\Omega)$ has an isolated zero at the origin if and only if the vector fields $\nabla f$ and
$W(f,\Omega)$ have an isolated zero at the origin.
\end{lemma}
{\it Proof.} See Lemma 7.4 $\hfill \Box$

\begin{lemma}
The form $\Omega$ has an isolated zero at the origin on $ \{f \ge 0\}$ and $\{ f \le 0\}$ if and only if
$Y(f,\Omega)$ has an isolated zero at the origin.
\end{lemma}
{\it Proof.} See Lemma 7.5. $\hfill \Box$

We can state the version of Theorem 7.6.
\begin{theorem}
Assume that $Y(f,\Omega)$ has an isolated zero at the origin. If $n$ is even, we have :
$$\displaylines{
\hbox{\em Ind}_{Rad}(\Omega,0,\{f\ge 0\})=\frac{1}{2} \Big[ \hbox{\em Ind}_{PH}(V(\Omega),0,\mathbb{R}^n) -
\hbox{\em Ind}_{PH}(W(f,\Omega),0,\mathbb{R}^n)+  \hfill \cr
\hfill \hbox{\em Ind}_{PH}(\nabla f,0,\mathbb{R}^n)-\hbox{\em Ind}_{PH}(Y(f,\Omega),0,\mathbb{R}^n) \Big], 
\quad \cr}$$
$$\displaylines{
\hbox{\em Ind}_{Rad}(\Omega,0,\{f\le 0\})=\frac{1}{2} \Big[ \hbox{\em Ind}_{PH}(V(\Omega),0,\mathbb{R}^n) -
\hbox{\em Ind}_{PH}(W(f,\Omega),0,\mathbb{R}^n)- \hfill \cr
\hfill \hbox{\em Ind}_{PH}(\nabla f,0,\mathbb{R}^n)+\hbox{\em Ind}_{PH}(Y(f,\Omega),0,\mathbb{R}^n) \Big].
\quad \cr
}$$
If $n$ is odd, we have :
$$\displaylines{
 \hbox{\em Ind}_{Rad}(\Omega,0,\{f\ge 0\})=\frac{1}{2} \Big[ \hbox{\em Ind}_{PH}(Y(f,\Omega),0,\mathbb{R}^n) +
\hbox{\em Ind}_{PH}(W(f,\Omega),0,\mathbb{R}^n)-  \hfill \cr 
\hfill \hbox{\em Ind}_{PH}(\nabla f,0,\mathbb{R}^n)\Big],
\quad \cr
}$$
$$\displaylines{
 \hbox{\em Ind}_{Rad}(\Omega,0,\{f\le 0\})=\frac{1}{2} \Big[ \hbox{\em Ind}_{PH}(Y(f,\Omega),0,\mathbb{R}^n) -
\hbox{\em Ind}_{PH}(W(f,\Omega),0,\mathbb{R}^n)+ \hfill \cr
\hfill   \hbox{\em Ind}_{PH}(\nabla f,0,\mathbb{R}^n)\Big].
}$$

\end{theorem}
{\it Proof.} Perturbing $f$ and $\Omega$ as in the previous theorems and using the same notations as in Theorem 7.6, we
find that :
$$\displaylines{
\quad \hbox{Ind}_{Rad}(\Omega,0,\{f \ge 0 \})+\hbox{Ind}_{Rad}(\Omega,0,\{f \le 0\})= \hfill \cr
\hfill \sum_j \hbox{Ind}_{PH}(\Omega,Q_j,f^{-1}(\delta))+\hbox{Ind}_{PH}(V(\Omega),0,\mathbb{R}^n) + \hfill \cr
\hfill 1-\chi(f^{-1}(\delta)\cap B_r), \quad \cr
}$$
and,
$$\displaylines{
\quad \hbox{Ind}_{Rad}(\Omega,0,\{f \ge 0 \})-\hbox{Ind}_{Rad}(\Omega,0,\{f \le 0\})= \hfill \cr
\hfill \sum_j  \hbox{sign}(\langle \nabla f (Q_j),V(\Omega)(Q_j)
\rangle) \hbox{Ind}_{PH}(\Omega,Q_j,f^{-1}(\delta))+ \qquad \qquad \cr
\hfill \hbox{sign}(-\delta) \hbox{Ind}_{PH}(V(\Omega),0,\mathbb{R}^n) - 
\left[ \chi\left( \{f\ge \delta\} \cap B_r^n\right) - \chi\left( \{ f \le \delta \} \cap B_r^n \right) \right]. \quad \cr
}$$
If $n$ is even, $\sum_j \hbox{Ind}_{PH} (\Omega,Q_j,f^{-1}(\delta))=-\hbox{Ind}_{PH}(W(f,\Omega),0,\mathbb{R}^n)$ and
$$1-\chi(f^{-1}(\delta)\cap B_r^n)=\hbox{Ind}_{PH}(\nabla f,0,\mathbb{R}^n),$$ and so :
$$\displaylines{
\quad \hbox{Ind}_{Rad}(\Omega,0,\{f \ge 0 \})+\hbox{Ind}_{Rad}(\Omega,0,\{f \le 0\})= \hfill \cr
\hfill -\hbox{Ind}_{PH}(W(f,\Omega),0,\mathbb{R}^n)+\hbox{Ind}_{PH}(V(\Omega),0,\mathbb{R}^n) + 
\hbox{Ind}_{PH}(\nabla f,0,\mathbb{R}^n). \quad \cr
}$$
Furthermore :
$$\displaylines{
\quad \hbox{Ind}_{PH}(Y(f,\Omega),0,\mathbb{R}^n)= \hfill \cr
\quad \quad -\sum_j  \hbox{sign}(\langle \nabla f (Q_j),V(\Omega)(Q_j)
\rangle) \hbox{Ind}_{PH}(\Omega,Q_j,f^{-1}(\delta))+\hfill \cr
\hfill \hbox{sign}(-\delta)\hbox{Ind}_{PH} (\Gamma(f,\Omega),0,\mathbb{R}^n) = \quad \cr
\quad  -\sum_j  \hbox{sign}(\langle \nabla f (Q_j),V(\Omega)(Q_j)
\rangle) \hbox{Ind}_{PH}(\Omega,Q_j,f^{-1}(\delta)) + \quad \quad \cr
\hfill \hbox{sign}(-\delta) \hbox{Ind}_{PH}(\nabla f,0,\mathbb{R}^n) -\hbox{sign}(-\delta)
\hbox{Ind}_{PH}(V(\Omega),0,\mathbb{R}^n). \quad \cr
}$$
Therefore :
$$\displaylines{
\quad \hbox{Ind}_{Rad}(\Omega,0,\{f \ge 0 \})-\hbox{Ind}_{Rad}(\Omega,0,\{f \le 0\})= \hfill \cr
\quad \quad  -\hbox{Ind}_{PH}(Y(f,\Omega),0,\mathbb{R}^n)-\hbox{sign}(\delta)
\hbox{Ind}_{PH}(\nabla f,0,\mathbb{R}^n) + \hfill \cr
\hfill \hbox{sign}(\delta) \hbox{Ind}_{PH}(V(\Omega),0,\mathbb{R}^n) +\quad \cr
\quad \quad  \hbox{sign}(-\delta) \hbox{Ind}_{PH}(V(\Omega),0,\mathbb{R}^n) -\big[\hbox{sign}(-\delta) \hbox{Ind}_{PH}
(\nabla f,0,\mathbb{R}^n) \big]= \hfill \cr
\hfill -\hbox{Ind}_{PH} (Y(f,\Omega),0,\mathbb{R}^n). \quad \cr
}$$
If $n$ is odd :
$$\displaylines{
\qquad \sum_j \hbox{Ind}_{PH}(\Omega,Q_j,f^{-1}(\delta))=  \hbox{Ind}_{Rad}(Y(f,\Omega),0,\mathbb{R}^n)+ \hfill \cr
\hfill \hbox{sign}(\delta)
\hbox{Ind}_{PH}(\nabla f,0,\mathbb{R}^n), \qquad \cr
}$$
$$1-\chi(f^{-1}(\delta)\cap B_r^n)=-\hbox{sign}(\delta) \hbox{Ind}_{PH}(\nabla f,0,\mathbb{R}^n),$$ and :
$$\hbox{Ind}_{Rad}(\Omega,0,\{f\ge 0\})+\hbox{Ind}_{Rad}(\Omega,0,\{f \le 0 \})=\hbox{Ind}_{PH}(Y(f,\Omega),0,\mathbb{R}^n).$$
Furthermore : $$\hbox{Ind}_{PH}(W(f,\Omega),0,\mathbb{R}^n)=\sum_j \langle \nabla f(Q_j), V(\Omega) (Q_j)\rangle 
\hbox{Ind}_{PH}(\Omega,Q_j,f^{-1}(\delta)),$$ so we obtain :
$$\displaylines{
\qquad \hbox{Ind}_{Rad}(\Omega,0,\{f \ge 0 \})- \hbox{Ind}_{Rad}(\Omega,0,\{ f \le  0\})= \hfill \cr
\hfill \hbox{Ind}_{PH}
(W(f,\Omega),0,\mathbb{R}^n) -\hbox{Ind}_{PH}(\nabla f,0,\mathbb{R}^n). \qquad \cr
}$$ $\hfill \Box$

\begin{corollary}
Let $g :(\mathbb{R}^n,0) \rightarrow (\mathbb{R},0)$ be an analytic function defined in the neighborhood of the origin with
$g(0)=0$. Let us assume that $\nabla g$ satisfies Condition $(P')$ and that $Y(f,dg)$ has an isolated zero at the
origin. If $n$ is even, we have :
$$\displaylines{
\quad \chi(\{f \ge 0\} \cap g^{-1}(\delta) \cap B_r^n)=1-\frac{1}{2}\Big[ \hbox{\em Ind}_{PH}(\nabla f,0,\mathbb{R}^n) +
\hbox{\em Ind}_{PH}(V(\Omega),0,\mathbb{R}^n) - \hfill \cr
\hfill \hbox{\em Ind}_{PH}(Y(f,dg),0,\mathbb{R}^n) \Big]- 
\frac{1}{2}\hbox{\em sign}(\delta) \hbox{\em Ind}_{PH}(W(f,dg),0,\mathbb{R}^n) , \quad \cr
\quad \chi(\{f \le 0\} \cap g^{-1}(\delta) \cap B_r^n)=1-\frac{1}{2}\Big[ \hbox{\em Ind}_{PH}(\nabla f,0,\mathbb{R}^n) +
\hbox{\em Ind}_{PH}(V(\Omega),0,\mathbb{R}^n) +  \hfill \cr
\hfill  \hbox{\em Ind}_{PH}(Y(f,dg),0,\mathbb{R}^n) \Big]- 
\frac{1}{2}\hbox{\em sign}(\delta) \hbox{\em Ind}_{PH}(W(f,dg),0,\mathbb{R}^n) . \quad \cr
}$$
If $n$ is odd, we have :
$$\displaylines{
\quad \chi(\{f \ge 0\} \cap g^{-1}(\delta) \cap B_r^n)=1-\frac{1}{2}\Big[- \hbox{\em Ind}_{PH}(\nabla f,0,\mathbb{R}^n)
+\ \hfill \cr
\hfill \hbox{\em Ind}_{PH}(Y(f,dg),0,\mathbb{R}^n) \Big] 
+\frac{1}{2}\hbox{\em sign}(\delta) \hbox{\em Ind}_{PH}(W(f,dg),0,\mathbb{R}^n) , \quad \cr
\quad \chi(\{f \le 0\} \cap g^{-1}(\delta) \cap B_r^n)=1-\frac{1}{2}\Big[ \hbox{\em Ind}_{PH}(\nabla f,0,\mathbb{R}^n) +
 \hfill \cr
\hfill \hbox{\em Ind}_{PH}(Y(f,dg),0,\mathbb{R}^n) \Big]
-\frac{1}{2}\hbox{\em sign}(\delta) \hbox{\em Ind}_{PH}(W(f,dg),0,\mathbb{R}^n) . \quad \cr
}$$
\end{corollary}

\noindent{\bf Examples}

$ \bullet $ In $\mathbb{R}^2$, let $f(x_1,x_2)=\frac{1}{2}(x_1^2-x_2^2)$ and $\Omega(x_1,x_2)=(x_1-x_2)dx_1+x_1dx_2$. 
It is easy to see that $\hbox{Ind}_{PH}(\nabla f,0,\mathbb{R}^n)=-1$ and $\hbox{Ind}_{PH}(V(\Omega),0,\mathbb{R}^n)=1$.
Moreover the computer gives that :$$\hbox{Ind}_{PH}(W(f,\Omega),0,\mathbb{R}^n)=2 \hbox{ and }
\hbox{Ind}_{PH}(Y(f,\Omega),0,\mathbb{R}^n)=0.$$
\noindent Applying Theorem 7.2 and Theorem 7.6, we obtain :
$$\hbox{Ind}_{Rad}(\Omega,0,f^{-1}(0))=1, 
\hbox{Ind}_{Rad}(\Omega,0,\{f \ge 0\})=1,$$
$$ \hbox{Ind}_{Rad}(\Omega,0,\{ f \le 0\})=1.$$

$ \bullet $ In $\mathbb{R}^2$, let $f(x_1,x_2)=x_1^3-x_2^2$ and $\Omega(x_1,x_2)=(x_1-x_2)dx_1+x_1dx_2$. 
It is easy to see that $\hbox{Ind}_{PH}(\nabla f,0,\mathbb{R}^n)=0$ and $\hbox{Ind}_{PH}(V(\Omega),0,\mathbb{R}^n)=1$.
Moreover the computer gives that : $$\hbox{Ind}_{PH}(W(f,\Omega),0,\mathbb{R}^n)=1 \hbox{ and }
\hbox{Ind}_{PH}(Y(f,\Omega),0,\mathbb{R}^n)=0.$$
\noindent Applying Theorem 7.2 and Theorem 7.6, we obtain :
$$\hbox{Ind}_{Rad}(\Omega,0,f^{-1}(0))=1, 
\hbox{Ind}_{Rad}(\Omega,0,\{f \ge 0\})=1,$$
$$ \hbox{Ind}_{Rad}(\Omega,0,\{ f \le 0\})=1.$$

$ \bullet $ In $\mathbb{R}^4$, let $f(x_1,x_2,x_3,x_4)=\frac{1}{2}(x_1^2-x_2^2+x_3^2+x_4^2)$ and $\Omega(x_1,x_2,x_3,x_4)=
x_4dx_1-x_1dx_2+x_2dx_3+x_3dx_4$. 
It is easy to see that $\hbox{Ind}_{PH}(\nabla f,0,\mathbb{R}^n)=-1$ and $\hbox{Ind}_{PH}(V(\Omega),0,\mathbb{R}^n)=1$.
The vector fields $W(f,\Omega)$ and $Y(f,\Omega)$ are given by :
$$\displaylines{
\qquad W(f,\Omega)(x_1,x_2,x_3,x_4)= 
   (\frac{1}{2}(x_1^2-x_2^2+x_3^2+x_4^2) , -x_1^2+x_3^2 , \hfill \cr
\hfill   -x_3x_4-x_1x_4+x_1x_2+x_2x_3, -x_4^2+2x_1x_3-x_2^2 ), \qquad \cr
}$$
and :
$$\displaylines{
\qquad Y(f,\Omega)(x_1,x_2,x_3,x_4)= 
   (\frac{1}{2}(x_1^2-x_2^2+x_3^2+x_4^2)(x_1x_4+x_1x_2+x_2x_3+x_3x_4) ,  \hfill \cr
\hfill  -x_1^2+x_3^2 , -x_3x_4-x_1x_4+x_1x_2+x_2x_3, -x_4^2+2x_1x_3-x_2^2 ), \qquad \cr
}$$
It is easy to check that these two mappings have an isolated zero at the origin in $\mathbb{R}^4$. This is not true any more in
$\mathbb{C}^4$ because the line in $\mathbb{C}^4$ through $(0,0,0,0)$ and
$(\frac{1}{\sqrt{2}},0,-\frac{1}{\sqrt{2}},i)$ is included in $W(f,\Omega)^{-1}(0)$.  Hence we can not use the program to
compute the indices of $W(f,\Omega)$ and $Y(f,\Omega)$. Nevertheless it is possible to compute them by hands. 
Since for $\varepsilon >0$ the point $(0,0,0,\varepsilon)$ has no preimage by $W(f,\Omega)$,
$\hbox{Ind}_{PH}(W(f,\Omega),0,\mathbb{R}^n)=0$. By $Y(f,\Omega)$, it has exactly two preimages : $(\alpha,0,\alpha,0)$ and
$-(\alpha,0,\alpha,0)$ where $\alpha=\sqrt{\frac{\varepsilon}{2}}$. At each of these point points, the jacobian determinant of
$Y(f,\Omega)$ is strictly positive. We conclude that $\hbox{Ind}_{PH}(Y(f,\Omega),0,\mathbb{R}^n)=2$. Applying Theorem 7.2 and Theorem 7.6, we obtain :
$$\hbox{Ind}_{Rad}(\Omega,0,f^{-1}(0))=-1, 
\hbox{Ind}_{Rad}(\Omega,0,\{f \ge 0\})=1,$$
$$ \hbox{Ind}_{Rad}(\Omega,0,\{ f \le 0\})=-1.$$

\section{Radial index  on semi-analytic curves}
In this section, we explain briefly how to compute the radial index of a 1-form on a semi-analytic curve defined as
the set of points on a $1$-dimensional complete intersection where some analytic inequalities are satisfied. 

First we give a characterization of the radial index on a subanalytic curve. Let $\mathcal{C}\subset \mathbb{R}^n$ be 
a subanalytic curve and let us assume that $0$ belongs to $\mathcal{C}$. Let $\Omega$ be a 1-form on $\mathbb{R}^n$ 
such that $0$ is an isolated zero of $\Omega$ on $\mathcal{C}$. Thus $\Omega$ defines an orientation on each half-branch of
 $\mathcal{C}\setminus \{0\}$. We say that a half-branch is inbound (resp. outbound) if the orientation is towards
(resp. away) from $0$. 
\begin{lemma}
If $\Omega$ has an isolated zero at $0$ on $\mathcal{C}$ then :
$$\hbox{\em Ind}_{Rad} (\Omega,0,\mathcal{C})=1-\# \{\hbox{\em inbound half-branches}\}.$$
\end{lemma}
{\it Proof.} Let $\tilde{\Omega}$ be a small perturbation of $\Omega$ which satisfies the three conditions stated before
Definition 4.2. Let $0<r'<r \ll 1$ be such that $\tilde{\Omega}$ is radial in $B_{r'}^n$ and coincides with $\Omega$ in the
neighborhood of $S_r^{n-1}$. Applying the Poincar\'e-Hopf theorem and the definition of the radial index and denoting by
$b(\mathcal{C})$ the number of half-branches of $\mathcal{C}\setminus \{0\}$ , we obtain :
$$b(\mathcal{C})=\hbox{Ind}_{Rad}(\Omega,0,\mathcal{C})-1 +b(\mathcal{C})+\# \{\hbox{inbound half-branches}\}.$$ $\hfill \Box$

Let $F=(f_1,\ldots,f_{n-1}) : (\mathbb{R}^n,0) \rightarrow (\mathbb{R}^{n-1},0)$ be an analytic mapping defined in the
neighborhood of the origin such that $F(0)=0$ and $0$ is isolated in 
$\{x \in \mathbb{R}^n \ \vert \ F(x)=0 \hbox{ and } \hbox{rank}[DF(x)]< n-1 \}$. This implies that $F^{-1}(0)$ is a curve
with an isolated singularity at the origin. Let $\Omega=a_1dx_1+\cdots+a_ndx_n$ be a smooth 1-form. 
Let $g : (\mathbb{R}^n,0) \rightarrow (\mathbb{R},0)$ be an analytic function defined in the neighborhood of the origin
such that $g(0)=0$. 
Let $V(\Omega)$ and $W(\Omega,g)$ be the following vector fields :
$$V(\Omega)=M(\Omega) \frac{\partial }{\partial x_1}+f_1\frac{\partial}{\partial x_2}+\cdots+
f_{n-1}\frac{\partial}{\partial x_{n}} ,$$
$$W(\Omega,g)=M(\Omega)g \frac{\partial }{\partial x_1}+f_1\frac{\partial}{\partial x_2}+\cdots+
f_{n-1}\frac{\partial}{\partial x_{n}},$$
where :
$$M(\Omega)=\left\vert \begin{array}{ccc}
a_1 & \ldots & a_n \cr
\frac{\partial f_1}{\partial x_1} & \ldots & \frac{\partial f_1}{\partial x_n} \cr
\vdots & \ddots & \vdots \cr
\frac{\partial f_{n-1}}{\partial x_1} & \ldots & \frac{\partial f_{n-1}}{\partial x_n} \cr
\end{array} \right\vert.$$

\begin{lemma}
The form $\Omega$ has an isolated
zero at $0$ on $F^{-1}(0)$ if and only if the vector field $V(\Omega)$ has an isolated zero at the origin.
\end{lemma}
{\it Proof.} It is clear. $\hfill \Box$

\begin{lemma}
The vector field $W(\Omega,g)$ has an isolated zero at the origin if and only if $\Omega$ has an isolated
zero at $0$ on $F^{-1}(0)$ and $g$ does not vanish on $F^{-1}(0) \setminus \{0\}$ in a neighborhood of the origin.
\end{lemma}
{\it Proof.} This is clear because $W(\Omega,g)$ has an isolated zero at the origin if and only if $V(\Omega)$ has an
isolated zero at the origin and $0$ is isolated in $g^{-1}(0) \cap F^{-1}(0)$.  $\hfill \Box$

Now let $V(dg)$ and $I$ be the following vector fields :
$$V(dg)=\frac{\partial(g,f_1,\ldots,f_{n-1})}{\partial(x_1,\ldots,x_n)} \frac{\partial }{\partial x_1}+f_1\frac{\partial}{\partial x_2}+\cdots+
f_{n-1}\frac{\partial}{\partial x_{n}} ,$$
$$I=\frac{\partial(\rho,f_1,\ldots,f_{n-1})}{\partial(x_1,\ldots,x_n)} \frac{\partial}{\partial x_1}+
f_1\frac{\partial}{\partial x_2}+\cdots+f_{n-1}\frac{\partial}{\partial x_{n}} ,$$
where $\rho(x)=x_1^2+\cdots+x_n^2.$
Note that $I=V(2\sum_i x_idx_i)=V(d\rho)$. 

\begin{lemma}
The vector field $I$ has an isolated zero at the origin.
Furthermore if $W(\Omega,g)$ has an isolated zero at the origin then $V(dg)$ has an isolated zero at the origin.
\end{lemma}
{\it Proof.} The first assertion is proved in [Sz1], Lemma 2.3. If $W(\Omega,g)$ has an isolated zero at the origin, then 
$0$ is isolated in $g^{-1}(0) \cap F^{-1}(0)$ by the previous lemma. We just have to apply Lemma 2.3 in [Sz1]. $\hfill
\Box$


\begin{theorem}
Assume that $W(\Omega,g)$ has an isolated zero at the origin. Then we have :
$$\displaylines{
\qquad \hbox{\em Ind}_{Rad}(\Omega,0,F^{-1}(0)\cap \{g \ge 0\})= 1+\frac{1}{2} \Big[
 \hbox{\em Ind}_{PH}(W(\Omega,g),0,\mathbb{R}^n) + \cr
\hfill  \hbox{\em Ind}_{PH}(V(\Omega),0,\mathbb{R}^n) 
 -\hbox{\em Ind}_{PH}(V(dg),0,\mathbb{R}^n)-\hbox{\em Ind}_{PH}(I,0,\mathbb{R}^n) \Big] .
\qquad \cr
}$$
$$\displaylines{
\qquad \hbox{\em Ind}_{Rad}(\Omega,0,F^{-1}(0)\cap \{g \le 0\})= 1+\frac{1}{2} \Big[
 -\hbox{\em Ind}_{PH}(W(\Omega,g),0,\mathbb{R}^n) + \cr 
 \hfill  \hbox{\em Ind}_{PH}(V(\Omega),0,\mathbb{R}^n) 
 +\hbox{\em Ind}_{PH}(V(dg),0,\mathbb{R}^n)-\hbox{\em Ind}_{PH}(I,0,\mathbb{R}^n) \Big] .
\qquad \cr
}$$

\end{theorem}
{\it Proof.} The proof of this theorem is very similar to the proofs of the theorems of the previous section so
we will not give all the details.

Let $(\delta,\alpha)$ be a regular value of $(F,g)$ such that $0\le \vert \alpha \vert \ll \vert \delta \vert \ll r$. 
We perturb $\Omega$ into $\tilde{\Omega}$ such that $\tilde{\Omega}$ is correct and non-degenerate on 
$F^{-1}(\delta) \cap B_r^n$, $F^{-1}(\delta) \cap \{g \ge \alpha \} \cap B_r^n$ and 
$F^{-1}(\delta) \cap \{g \le \alpha \} \cap B_r^n$. Denoting by $Q_1,\ldots,Q_s$ the singular points of
$\tilde{\Omega}$ on $F^{-1}(\delta)$ lying in $\mathring{B_r^n}$ and using the Poincar\'e-Hopf theorem, we find that :
$$\displaylines{
\quad \chi(F^{-1}(\delta) \cap B_r^n)=\sum_{i=1}^s \hbox{Ind}_{PH}(\tilde{\Omega},Q_i,F^{-1}(\delta))+ \hfill \cr
\hfill 2-\hbox{Ind}_{Rad}(\Omega,0,F^{-1}(0)\cap \{g \ge 0 \})-\hbox{Ind}_{Rad}(\Omega,0,F^{-1}(0)\cap \{ g \le 0 \}).
\quad \cr
}$$
By Lemma 2.1, it is easy to see that :
$$\sum_{i=1}^s \hbox{Ind}_{PH}(\tilde{\Omega},Q_i,F^{-1}(\delta))=\hbox{Ind}_{PH}(V(\Omega),0,\mathbb{R}^n).$$ 
Furthermore, $\chi(F^{-1}(\delta) \cap
B_r^n)=\hbox{Ind}_{PH}(I,0,\mathbb{R}^n)$ (see [AFS], [AFN], [Sz1]). Hence:
$$\displaylines{
\quad \hbox{Ind}_{Rad}(\Omega,0,F^{-1}(0)\cap \{g \ge 0 \})+\hbox{Ind}_{Rad}(\Omega,0,F^{-1}(0)\cap \{ g \le 0 \})=
2+ \hfill \cr
\hfill \hbox{Ind}_{PH}(V(\Omega),0,\mathbb{R}^n)-\hbox{Ind}_{PH}(I,0,\mathbb{R}^n). \quad \cr
}$$
Let us write $F^{-1}(\delta) \cap g^{-1}(\alpha) \cap B_r^n =\{P_1,\ldots,P_r\}$. Since $\tilde{\Omega}$ is correct on 
$F^{-1}(\delta) \cap \{g \ge \alpha \} \cap B_r^n$ and 
$F^{-1}(\delta) \cap \{g \le \alpha \} \cap B_r^n$, for each $j \in \{1,\ldots,r\}$ there exists $\lambda_j \not= 0$ such that
: 
$$\Omega_{\vert F^{-1}(\delta)} (P_j)=\lambda_j dg_{\vert F^{-1}(\delta)}(P_j).$$
By the Poincar\'e-Hopf theorem for manifolds with corners, we have :
$$\displaylines{
 \chi(F^{-1}(\delta) \cap \{g \ge \alpha\} \cap B_r^n)=
\sum_{i \ \vert \ g(p_i)> \alpha} \hbox{Ind}_{PH}(\tilde{\Omega},Q_i,F^{-1}(\delta))
+\# \{j \ \vert \ \lambda_j >0\} + \hfill \cr
\hfill 1-\hbox{Ind}_{Rad}(\Omega,0,F^{-1}(0)\cap \{ g \ge 0 \}), \quad \cr 
 \chi(F^{-1}(\delta) \cap \{g \le \alpha\} \cap B_r^n)=
\sum_{i \ \vert \ g(p_i)< \alpha} \hbox{Ind}_{PH}(\tilde{\Omega},Q_i,F^{-1}(\delta))
+\# \{j \ \vert \ \lambda_j <0\} + \hfill \cr
\hfill 1-\hbox{Ind}_{Rad}(\Omega,0,F^{-1}(0)\cap \{ g \le 0 \}). \quad \cr
}$$
This leads to :
$$\displaylines{
\hbox{Ind}_{Rad}(\Omega,0,F^{-1}(0)\cap \{ g \ge 0 \})-\hbox{Ind}_{Rad}(\Omega,0,F^{-1}(0)\cap \{ g \le 0 \})= 
\cr
\quad \quad  \sum_i \hbox{sign } (g(Q_i)-\alpha) \hbox{Ind}_{PH}(\tilde{\Omega},Q_i,F^{-1}(\delta)) + \sum_j \hbox{sign } \lambda_j
 \hfill \cr
\hfill -\Big[ \chi(F^{-1}(\delta) \cap \{g \ge \alpha\} \cap B_r^n)-\chi(F^{-1}(\delta) \cap \{g \le \alpha\} \cap B_r^n) \Big]. \quad \cr
}$$
A computation based on Cramer's rules shows that for each $j \in \{1,\ldots,r\}$ :
$$\hbox{sign } \lambda_j =\hbox{sign } \frac{\partial (g,f_1,\ldots,f_{n-1})}{\partial( x_1,\ldots,x_n)} (P_j)
\hbox{sign } M(\tilde{\Omega})(P_j).$$
Furthermore if $b_+(g)$ (resp. $b_-(g)$) is the number of half-branches of $F^{-1}(0)$ on which $g>0$ (resp. $g<0$), we
have by Theorem 3.1 in [Sz1] : 
$$\displaylines{
\quad \chi(F^{-1}(\delta) \cap \{g \ge \alpha\} \cap B_r^n)-\chi(F^{-1}(\delta) \cap \{g \le \alpha\} \cap B_r^n)= \hfill \cr
\hfill \frac{1}{2}
\left( b_+(g)-b_-(g) \right)= \hbox{Ind}_{PH} (V(dg),0,\mathbb{R}^n). \quad \cr
}$$
Collecting all these informations, we obtain :
$$\displaylines{
\quad \hbox{Ind}_{Rad}(\Omega,0,F^{-1}(0)\cap \{ g \ge 0 \})-\hbox{Ind}_{Rad}(\Omega,0,F^{-1}(0)\cap \{ g \le 0 \})= \hfill
\cr 
\hfill \hbox{Ind}_{PH}(W(\Omega,g),0,\mathbb{R}^n) -\hbox{Ind}_{PH}(V(dg),0,\mathbb{R}^n). \quad \cr
}$$
$\hfill \Box$

Let us apply this theorem when $g=\rho$. In this case, $V(dg)=V(d\rho)=I$ and $W(\Omega,\rho)$ has an
isolated zero at the origin if and only if $V(\Omega)$ has an isolated zero at the origin. Futhermore, in this situation,
these two vector fields have the same index at the origin because on a small sphere they never point in opposite
directions. We can state :
\begin{corollary}
Assume that $V(\Omega)$ has an isolated zero at the origin. Then $\Omega$ has an isolated zero at the origin on
$F^{-1}(0)$ and :
$$\hbox{\em Ind}_{Rad}(\Omega,0,F^{-1}(0))=1 +\hbox{\em Ind}_{PH}(V(\Omega),0,\mathbb{R}^n)-\hbox{\em
Ind}_{PH}(I,0,\mathbb{R}^n).$$ $\hfill \Box$
\end{corollary}

We will generalize Theorem 8.5 to the case of a closed semi-analytic curve defined by several sign conditions. 
More precisely let $g_1,\ldots,g_k :(\mathbb{R}^n,0) \rightarrow (\mathbb{R},0)$ be analytic functions defined in 
the neighborhood of the origin such that $g_j(0)=0$, for $j \in \{1,\ldots,k\}$. For each
$\alpha=(\alpha_1,\ldots,\alpha_k) \in \{0,1\}^k$, let us define the vector fields $W(\Omega,\alpha)$ 
and $V(\alpha)$ in the following way :
$$\hbox{if } \alpha \not= (0,\ldots,0), \ W(\Omega,\alpha)=W(\Omega,g_1^{\alpha_1}\cdots g_k^{\alpha_k})\ , \ 
V(\alpha)=V(d(g_1^{\alpha_1}\cdots g_k^{\alpha_k})),$$
$$W(\Omega,(0,\ldots,0))=V(\Omega) \ , \  V((0,\ldots,0))=I.$$
For each $\epsilon=(\epsilon_1,\ldots,\epsilon_k) \in \{0,1\}^k$, let $\mathcal{C}(\epsilon)$ be the semi-analytic curve
defined by :
$$\mathcal{C}(\epsilon)=F^{-1}(0) \cap \{(-1)^{\epsilon_1}g_1 \ge 0,\ldots,(-1)^{\epsilon_k}g_k \ge 0\}.$$

\begin{theorem}
If $W(\Omega,(1,\ldots,1))$ has an isolated zero at the origin then for every $\alpha \in \{0,1\}^k$, $W(\Omega,\alpha)$ and 
$V(\alpha)$ have an isolated zero at the origin. In this case, for every $\epsilon=(\epsilon_1,\ldots,\epsilon_k) \in
\{0,1\}^k$, we have :
$$\displaylines{
\qquad \hbox{\em Ind}_{Rad}(\Omega,0,\mathcal{C}(\epsilon))=1+ \hfill \cr
\hfill \frac{1}{2^k}
\sum_{\alpha \in \{0,1\}^k} (-1)^{\epsilon \cdot \alpha} \Big[ \hbox{\em Ind}_{PH} (W(\Omega,\alpha),0,\mathbb{R}^n)-
\hbox{\em Ind}_{PH} (V(\alpha),0,\mathbb{R}^n) \Big],\qquad \cr
}$$
where $\epsilon \cdot \alpha=\sum_{i=1}^k \epsilon_i \alpha_i $. 
\end{theorem}
{\it Proof.} The first affirmation is easy to check using Lemma 8.3 and Lemma 8.4. 
Let us prove the formula for $\hbox{Ind}_{Rad}(\Omega,0,\mathcal{C}(\epsilon))$ by induction on $k$. 
For $k=1$, this is Theorem 8.5. Now assume that $k>1$. Let us fix $\epsilon'=(\epsilon",\epsilon_{k-1}) \in \{0,1\}^{k-1}$ and let
$\epsilon^0=(\epsilon',0)$ and $\epsilon^1=(\epsilon',1)$. Since the radial index is $1-\# \{\hbox{ inbound half-branches}    
 \}$, we get
:
$$\hbox{Ind}_{Rad}(\Omega,0,\mathcal{C}(\epsilon^0))+\hbox{Ind}_{Rad}(\Omega,0,\mathcal{C}(\epsilon^1))=
1+\hbox{Ind}_{Rad}(\Omega,0,\mathcal{C}(\epsilon')),$$
$$\displaylines{
\qquad \hbox{Ind}_{Rad}(\Omega,0,\mathcal{C}(\epsilon^0))-\hbox{Ind}_{Rad}(\Omega,0,\mathcal{C}(\epsilon^1))= \hfill \cr
\qquad \qquad  \hbox{Ind}_{Rad}(\Omega,0,\mathcal{C}(\epsilon")\cap \{ (-1)^{\epsilon_{k-1}} g_{k-1}g_k \ge 0\})- \hfill \cr
\hfill \hbox{Ind}_{Rad}(\Omega,0,\mathcal{C}(\epsilon")\cap \{ (-1)^{1} g_k \ge 0\}). \qquad \cr
}$$
It is enough to use the inductive hypothesis to conclude. $\hfill \Box$

\noindent{\bf Example}

$ \bullet $ In $\mathbb{R}^3$, let $f_1(x_1,x_2,x_3)=x_1^2+x_2^2-x_3^2$, $f_2(x_1,x_2,x_3)=x_1x_2$ and 
$g(x_1,x_2,x_3)=x_1^2-3x_2^2+x_3^2$. Let $\Omega(x_1,x_2,x_3)=(x_3^2+x_2)dx_1+x_1dx_2+(x_3^2-x_2^2)dx_3$. 
The computer gives that : $$\hbox{Ind}_{PH}(I,0,\mathbb{R}^n)=4, 
\hbox{Ind}_{PH}(V(dg),0,\mathbb{R}^n)=0,$$
$$\hbox{Ind}_{PH}(V(\Omega),0,\mathbb{R}^n)=\hbox{Ind}_{PH}(W(\Omega,g),0,\mathbb{R}^n)=0.$$
\noindent Applying Theorem 8.5 and Corollary 8.6, we obtain :
$$\hbox{Ind}_{Rad}(\Omega,0,F^{-1}(0)\cap \{g \ge 0\})=
\hbox{Ind}_{Rad}(\Omega,0,F^{-1}(0)\cap \{g \le 0\})=-1,$$ and :
$$\hbox{Ind}_{Rad}(\Omega,0,F^{-1}(0))=-3.$$

Let us end this section with a remark on a paper of Montaldi and van Straten.
In [MvS], Montaldi and van Straten study 1-forms on singular curves. They first consider the case of a meromorphic form
$\alpha$ on a reduced analytic curve $\mathcal{C}$ with base point $p$. They say that $\alpha$ is a finite form if its
restriction to each branch is not identically zero. When $\alpha$ is finite the define two ``ramification modules" which
are finite dimensional vector spaces and they prove that the difference of their dimensions is preserved under deformation
of the form and the curve. Then they consider the real case. They say that a real analytic curve $\mathcal{C}$ with base 
point $p$ is reduced if its complexification is and that a 1-form on $\mathcal{C}$ is meromorphic and finite if its
complexification is. In Theorem 2.1
and Corollary 2.2, they give formulas which express the number of outbound half-branches at $p$ and the number of inbound
half-branches at $p$ in terms of signatures of non-degenerate quadratic forms defined on appropriate vector spaces. Therefore,
by Lemma 8.1, Montaldi and van Straten's
results provide an Eisenbud-Levine type formula for the radial index of a meromorphic 1-form on a real reduced analytic
curve-germ.

\end{document}